\documentclass[12pt,a4paper,leqno]{article}
\usepackage{amscd,amsmath,amssymb,amsfonts}
\input xy
\xyoption{all}

\newtheorem{theo}{Theorem}[section]
\newtheorem{exam}[theo]{Examples}
\newtheorem{lem}[theo]{Lemma}

\newtheorem{prop}[theo]{Proposition}

\newtheorem{coro}[theo]{Corollary}

\newtheorem{rem}[theo]{Remark}

\newtheorem{defi}[theo]{Definition}

\def\noi{\noindent}

\def\dim{\mathrm{dim}}

\def\zn{{\mathbb Z/n\mathbb Z}}

\def\n{{\mathbb N}}
\def\z{{\mathbb Z}}

\def\r{{\mathbb R}}
\def\c{{\mathbb C}}

\def\ker{\mathrm{ker}}
\def\coker{\mathrm{coker}}

\def\cCH{\mathcal{CH}}

\def\cH{{\mathcal H}}
\def\cK{{\mathcal K}}

\def\cS{{\mathcal S}}

\def\cM{{\mathcal{M}}}

\def\Tl{{Tor^\ell}}
\def\Cl{{Cotor^\ell}}

\def\v{\vspace{0,5cm}}

\begin{document}

\title{\Large Rigidity theorems for $\mathcal{K}$- and $\mathcal{H}$-cohomology and other functors}

\author{Uwe Jannsen}

\maketitle

%\begin{abstract}
\noi
Suslin \cite{Su} proved that for an extension $K/k$ of algebraically closed fields the induced maps
$K_m(k)[n] \rightarrow K_m(K)[n]$ and $K_m(k)/n \rightarrow K_m(K)/n$ for the higher algebraic
$K$-groups are isomorphisms, where $A[n]$  is the subgroup of $n$-torsion in an abelian group $A$,
and $A/n = A/nA$, by definition. In this paper we generalize this to other functors and other field extensions.
%\end{abstract}

%\tableofcontents

%\setcounter{section}{-1}

%\section{Introduction}

\v\noindent
More precisely, we will show rigidity results for the following functors on separated noetherian schemes $Y$ for a field $k$, which are contravariant for flat morphisms.

\v\noi
List of rigid functors

\v\noi (1) $V(Y) = H^i(Y_{\acute{e}t},\zn(j))$ (\'etale cohomology), where $i,j\in \z$, $n\in \n$ is invertible in $k$, and $\zn(j) = \mu^{\otimes j}_n$,
for the sheaf $\mu_n$ of $n$-th roots of unity on $V$.

\v\noi (2) $V(Y) = K_m(Y)$ (algebraic $K$-theory) for $m\in \n_0$.

\v\noi (3) $V(Y) = CH_{r,s}(Y)$ (Gillet's Chow homology \cite{Gi}) for $r,s\in\z$.

\v\noi (4) $V(Y) = CH^r(Y,s)$ (Bloch's higher Chow groups \cite{Bl}) for $r,s \in \n_0$, which contains the classical Chow groups $CH^r(Y)$
as the special case $CH^r(Y,0)$.

\v\noi (5) $V(Y) = H^i(Y,\cK_m)$ ($\cK$-cohomology) for $i,m \in \n_0$, the $i$-th Zariski cohomology of the Zariski sheaf associated to the
presheaf $U \mapsto K_m(U)$.

\v\noi (6) $V(Y) = H^i(Y,\cCH^r(s))$, where $\cCH^r(s)$ is the Zariski sheaf on $Y$ associated to the presheaf $U \mapsto CH^r(U,s)$.

\v\noi (7) $V(Y) = H^i(Y,\cH^j(\zn(k)))$ ($\cH$-cohomology), the $i$-th cohomology of the Zariski sheaf on $Y$ associated
to the presheaf $U \mapsto H^j(U_{\acute{e}t},\zn(k))$, for $i\in \n_0$ and $j, k \in \z$.

\v\noi (8) $V(Y) = E^2_{p,q}(Y,\zn)$, the $E^2$-terms of the niveau spectral sequence associated to the \'etale
cohomology in (1).

\v\noi (9) $V(Y) = H_a(C^{r,s}_n(Y))$, where $C^{r,s}_n(Y)$ is a homological complex of Galois cohomology groups introduced by Kato in \cite{Ka}.

\v On the other hand, we consider field extensions $K/K^0$ over $k$ with the following property,
where $n$ is a natural number.

\begin{defi}\label{Definition0.1}
A field extension $K/L$ is called $n$-rigid, if the following two properties hold
for any intermediate field $L$, $K^0 \subset L \subset K$, which is algebraically
closed in $K$, and any smooth curve $C$ over $L$.

\smallskip\noi
(a) If $C(K)\neq \emptyset$, then $C(L) \neq \emptyset$.

\smallskip\noi
(b) If $\overline{C}$ is the regular proper model of $C$ over $L$ and any closed,
not necessarily reduced closed subscheme $C_\infty \subset C$ the map
$$
J_{C_\infty}(C)(L)/n \longrightarrow J_{C_\infty}(C)(K)/n
$$
is surjective, where $J_{C_\infty}(C)$ is the Rosenlicht generalized Jacobian
(\cite{Ro}, \cite{Se} Chap.V and \cite{CC}),
which is a commutative, smooth, geometrically connected group variety over $L$.
\end{defi}

\smallskip\noi
\begin{exam}\label{Examples 0.2}
Examples for $n$-rigid field extensions are

\smallskip\noi (i) An extension $K/K^0$ of algebraically closed fields, where $n$ is arbitrary.

\smallskip\noi(ii) $K = \r$ and $K^0$ is dense and algebraically closed in $\r$, and $n$ is arbitrary.

\smallskip\noi(iii) $K$ is a complete discrete valuation field, $K^0$ is algebraically closed in $L$
and dense in $L$ for the valuation topology, and $n$ is invertible in $L$.
\end{exam}

\smallskip\noi
An application we have in mind for $(iii)$ is the situation where $K$ is a global field (i.e., a number field or a function field in one variable
over a finite field), and where $L = K_v$ is the completion of $K$ with respect to a discrete valuation $v$ of $K$, and $L^0$ is the
Henselization $K_{(v)}$ of $K$ at $v$. An example for $(ii)$ is the situation where $L$ is the completion of $K$ at a real place $v$
(hence isomorphic to the field $\r$ of real numbers), and $L^0$ is the associated real closure of $K$ in $L$, the algebraic elements in the extension $L/K$. An example for $(i)$ is the situation where $L$ is the completion of $K$ at a complex place (hence isomorphic to $\c$), and $L^0$ is
the algebraic closure of $K$ in $L$.

\smallskip\noi
In particular we show the following.

\begin{theo}\label{Theorem 0.3} If $V$ is one of the functors (1) - (9) in
	the list of rigid functors and the field extension $L^0 \subset L$ is one of the field
extensions (i) - (iii) in the examples 0.2, and $X$ is a separated variety over $L^0$, then, with $n$ as in the respective cases,
the restriction maps
$$
V(X)[n] \mathop{\longrightarrow}\limits^{\cong}  V(X_L)[n] \hspace{1cm} \mbox{    and    } \hspace{1cm} V(X)/n  \mathop{\longrightarrow}\limits^{\cong}  V(X_L)/n $$
are isomorphisms, where $X_L = X\times_{L^0}L$, $A[n] = ker (A \mathop{\longrightarrow}\limits^{n} A)$, and
$A/n = coker(A \mathop{\longrightarrow}\limits^{n} A)$.
\end{theo}

\smallskip\noi
\begin{rem}\label{rem1} (a) This is known for (1), where case (i) is well-known and implies (ii) and (iii) by Galois descent.

\v\noi
(b) For (2), case (i), and $Y = Spec(L^0)$ this is the result of Suslin quoted in the beginning.

\v\noi
(c) For (3) and case (i) this was proved by F. Lecomte \cite{Lec}.

\v\noi
(d) The cases (ii) and (iii) of (b) and (c) do not follow automatically - because these theories do not have Galois descent.

\v\noi
(e) The cases (ii) and (iii) of (9) were the motivation for this paper, and turned out to be the most
difficult ones.
\end{rem}

\noi
\begin{exam}\label{ex1} Let $K$ be a global field, let $v$ be a non-archimedean place.

\smallskip\noi
(a) For case (3) and $r = s$ we consider the Chow groups $H^r(X,\cK_r) = CH^r(X)$, i.e., the classical Chow groups
of a smooth projective variety $X$ over a field. In case (iii) we get isomorphisms, for a global field $K$, and
an integer $n$ invertible in $K$, and a non-archimedean place $v$ of $K$:
$$
CH^r(X_{K_{(v)}})[n] \mathop{\longrightarrow}\limits^{\cong}  CH^r(X_{K_v})[n] \hspace{1cm} \mbox{    and    }
\hspace{1cm} CH^r(X_{K_{(v)}})/n  \mathop{\longrightarrow}\limits^{\cong}  CH^r(X_{K_v})/n \,.
$$

\v\noi (b) Let $X$ be a smooth projective curve over $K$ and consider the ``residue map'' for $n$ invertible in $K$
$$
H^i(K(X),\zn(j)) \mathop{\longrightarrow}\limits^{\alpha} \mathop{\oplus}\limits_{x\in\mid X \mid} H^{i-1}(k(x),\zn(j))\,,
$$
where $\mid X \mid$ is the set of closed points of $X$. Then one has $\ker(\alpha) = H^0(X,\cH^i_n(j))$, and by the above
this group for $X_{K_{(v)}}$ is isomorphic to the one for $X_{K_v}$. The same holds for $\coker(\alpha) = H^1(X,\cH^i_n(j))$.

\v\noi (c) For a field $F$ one has a functorial isomorphism $K^M_m(F) \cong CH^m(F,m)$ between the Milnor $K$-group and the written
Bloch higher Chow group, by work of Nesterenko and Suslin \cite{NeSu}, see also Totaro \cite{To}. Hence one gets rigidity
for Milnor $K$-theory mod $n$ and the $n$-torsion of Milnor $K$-theory, for all three cases (i), (ii) and (iii).
\end{exam}

%The second main result of this paper is the following.
%
%\begin{theo}\label{Theorem02}
%Let $k$ be a field and let $R$ be the Henselization of the local ring of a smooth $k$-variety
%at a $k$-rational point. Then, for any smooth $k$-scheme $X$ and any integer $n$ invertible in $k$
%the restriction maps
%$$
%V(X\times_kR)[n] \mathop{\longrightarrow}  V(X)[n] \hspace{1cm} \mbox{    and    } \hspace{1cm} V(X\times_kR)/n  \mathop{\longrightarrow}\limits^{\cong}  V(X)/n %$$
%are isomorphisms, if $V$ is one of the functors (1) - (5).
%\end{theo}

%\begin{exam}

%\v\noi (a) This result is well-known for example (1) (\'etale cohomology with coefficients invertible in $k$),
%where it is known by the base change theorem in \'etale cohomology.
%
%\v\noi (b) For example (2) (algebraic $K$-theory) the result is known by results of Gabber, Gillet, and Thomason %\cite{GiTho}
%\end{exam}

\noi
The above results are implied by some more general theorems. We introduce a notion of rigid functors (see Definition 1.1)
and sufficiently rigid functors (see Definition 3.1), and prove that the above rigidity result hold for such functors,
justifying their names in retrospective. Then we show how to produce such functors from twisted Poincar\'e duality
theories as introduced by Bloch and Ogus \cite{BO}.

\section{Rigid Functors}

\smallskip\noindent Let $\mathcal S$ be a category of schemes, and let $\mathcal S^{flat}$ be the category just
endowed with the flat morphisms.

\begin{defi}\label{Definition1.1} A contravariant functor $V$
on ${\mathcal S}^{flat}$ with values in the category Ab of abelian groups
is called rigid, if it satisfies the following properties,
provided the occurring schemes and morphisms are in $\mathcal S$.

\v\noi (a) For any flat finite morphism $\pi: X\to Y$ there
is a transfer morphism $\pi_\ast: V(X)\to V(Y)$, such that for
another flat finite morphism $\rho: Y\to Z$ one has
$(\rho\pi)_\ast = \rho_\ast \pi_\ast$.

\v\noi (b) For every cartesian diagram of schemes
$$
\xymatrix{X' \ar[r]^{f'} \ar[d]_{\pi'} &  X\ar[d]^{\pi}& \\
Y'\ar[r]^f & Y & ,}
$$
with $\pi$ finite and flat, one has $f^\ast\pi_\ast = {\pi'}_\ast
f'^\ast:V(X)\to V(Y')$.

\v\noi (c) If $X=X_1\amalg X_2$, then the immersions
$\pi_i:X_i\hookrightarrow X
 \ (i=1,2)$ induce an isomorphism
$$
(\pi^\ast_1,\pi_2^\ast):V(X)\stackrel{\sim}{\longrightarrow}
V(X_1)\oplus V(X_2)
$$
with inverse $(\pi_1)_\ast + (\pi_2)_\ast$.

\v\noi (d) If $X_m=X\times_{\mathbb Z}$Spec $(\mathbb Z[T]/(T^m))$
is the $m$-fold thickening of $X$, then for the morphism
$\pi:X_m\to X$ one has $\pi^\ast\pi_\ast =$ multiplication by $m$.

\v\noi (e) If $\mathbb A^1_X=X\times_{\mathbb Z} {\mathbb{A}}^1_{\mathbb Z}$ is the affine line over $X$, then the projection
$p:\mathbb A^1_X\to X$ induces an isomorphism
$p^\ast:V(X)\stackrel{\sim}{\longrightarrow} V(\mathbb A^1_X)$.

\v\noi (f) Let $i\mapsto X_i$ be a filtered projective system of
schemes such that the transition morphisms $X_i\to X_j$ are
affine, and let
$X=\mathop{lim}\limits_{\longleftarrow} X_i$.
Then the canonical map
$$
\mathop{lim}\limits_{\longrightarrow} V(X_i)\longrightarrow V(X)
$$
is an isomorphism.
\end{defi}

\begin{theo}\label{Theorem 1.2} Let $K/K^0$ be an $n$-rigid field extension
(see Definition 0.2), and let $V$ be a rigid functor on the category
of all noetherian $K^0$-schemes, such that the value groups $V(Y)$ are $n$-torsion
groups. Then the morphisms
$$
V(X)[n] \mathop{\longrightarrow}\limits^{\cong}  V(X_L)[n] \hspace{1cm} \mbox{    and    } \hspace{1cm} V(X)/n  \mathop{\longrightarrow}\limits^{\cong}  V(X_L)/n $$
are isomorphisms.
\end{theo}

\bigskip\noindent  \textbf{Proof}  For the surjectivity it suffices to
prove

\medskip\noindent {\bf Claim 1:} \ If $F$ is a function field over
$K^0$, contained in $K$, then
$$
Im(V(F)\to V(K))\ \subseteq \ Im(V(K^0)\to V(K))\,.
$$

\medskip\noindent In fact, by the limit property 1.1 (f), $V(K)$ is
generated by the images of the maps $V(F)\to V(K)$ for all such
fields $F$.

\medskip\noindent We prove Claim 1 by induction on $d=deg.tr.(F/K^0)$, the
degree of transcendence of $F$ over $K^0$. If $d=0$, then
necessarily $F=K^0$, since $K^0$ is algebraically closed
in $K$, and the claim is trivially true. If $d>0$, then by Noether
normalization and separability of $F$ over $K^0$, there exists
a function field $F_1$ with $deg.tr.(F_1/K^0)=d-1$ and a
smooth, geometrically irreducible curve $C_1$ over $F_1$ such that
$F=F_1(C_1)$, the function field of $C_1$. Let $\tilde F_1$ be the
algebraic closure of $F_1$ in $K$ and let $\tilde C_1 =
C_1\times_{F_1}\tilde F_1$, and $\tilde F = \tilde F_1(\tilde
C_1)$. Then it suffices to show

\medskip\noindent {\bf Claim 2}: $Im(V(\tilde F)\to V(K)\subseteq
Im(V(\tilde F_1)\to V(K))$.

\medskip\noindent In fact, by 1.1 (f) every $\alpha\in Im(V(\tilde
F_1)\to V(K))$ lies in $Im(V(F_2)\to V(K))$ for a function field
$F_2$ over $K$ with $deg.tr.(F_2/K^0)=d-1$, so by induction
$\alpha$ lies in $Im(V(K^0)\to V(K))$.

\medskip\noindent In other words, it suffices to prove Claim 1 for a field
$L\subseteq K$ in place of $K^0$, which is
algebraically closed in $K$ but for which $K$ is not necessarily
separable over $L$, and for a function field $F$ with
$deg.tr.(F/L) = 1$ which is separable over $L$.

\medskip Let $\alpha\in Im(V(F)\to V(K))$. By 1.1 (f) there is a
smooth, geometrically irreducible curve $C$ over $L$ with function
field $L(C)=F$ such that $\alpha\in Im(V(C)\to V(K))$. Let
$Div(C)$ be the group of divisors on $C$, i.e., the free abelian
group on the closed points $x\in C$. Consider the bilinear pairing
$$
Div(C)\times V(C)\longrightarrow V(L)
$$
defined by sending $(x,\beta)$ to $(\pi_x)_\ast
\varphi_x^\ast(\beta)$, where $\varphi_x:Spec \ \kappa(x)\to C$
and $\pi_x: Spec \ \kappa(x)\to Spec \ L$ are the canonical
morphisms. Denote by $\overline C$ the regular proper model of $C$
and set $C_\infty = \overline C\setminus C$. Let $f$ be a
meromorphic function on $\overline C$ which is defined and equal
to one on $C_\infty$ (i.e., for $x\in C_\infty$, $f$ lies in
$\mathcal O_{\overline C,x}$ and its image in $\kappa(x)$ is 1).
Then the principal divisor (f) lies in the kernel of the above
pairing. In fact, $f$ defines a covering $\pi: C'\to \mathbb A^1_L
= \mathbb P^1_L-{1}$, where $C'$ is obtained from $C$ by deleting
the points where $f$ is defined and equal to one. Now one has a
commutative diagram of pairings
$$
\begin{array}{ccccc}
Div(C) & \times & V(C) & \rightarrow & V(L)\\
\cup\hspace{-1mm}\mid & & \downarrow j^\ast & & \parallel\\
Div(C') & \times & V(C') & \rightarrow & V(L)\\
\pi^\ast\uparrow & & \downarrow \pi_\ast & & \parallel\\
Div(\mathbb A^1_L) & \times & V(\mathbb A^1_L) & \rightarrow &
V(L)\; ,
\end{array}
$$
where $j^\ast$ is induced by the open immersion $j: C' \hookrightarrow C$.
Indeed, recall that for a closed point $x\in \mathbb A^1_L$ one has
$$
\pi^\ast(x)=\mathop{\sum}\limits_{\pi(y)=x}e(y/x)\cdot x \quad ,
$$
where $e(y/x) = \hbox{length}(\mathcal O_{C,y}\otimes\kappa(x))$
is the ramification index of $y$ over $x$ (the tensor product is
over $\mathcal O_{\mathbb A^1_L,x})$. Consider the following
cartesian diagram
$$
\xymatrix{ C'_x\ar[d]_{\pi '}\ar[r]^{\varphi '{_x}} & C'\ar[d]^{\pi} \\
Spec\,\kappa(x)\ar[r]^{\varphi_x} & {\mathbb A{^1}{_F}} & .}
$$
Then $\varphi^\ast_x\pi_\ast = (\pi')_\ast(\varphi'_x)^\ast$ by
1.1 (b), and it suffices to show that
$$
(\pi')_\ast = \mathop{\sum}\limits_{\pi(y)=x} \ e(y/x)(\pi'\alpha_y)_\ast
\alpha^\ast_y \quad ,
$$
where $\alpha_y:Spec \ \kappa(y)\to C'_x$ is the canonical morphism. But since
$C'$ is smooth over $L$, we have
$$
C'_x\cong\scriptstyle\coprod\limits_{\pi(y)=x}\textstyle
Spec(\kappa(y)[T]/(T^{e(y/x)})) \ ,
$$
so the wanted equality follows easily from 1.1 (c) and (d), cf.
also the proof of Lemma 1.8 below.

\medskip\noindent Now by definition $(f)
= \pi^\ast(0-\infty)$; hence the commutative diagram of pairings
shows that $\psi((f),\beta)\in V(L)$ coincides with
$\psi'(0-\infty,\pi_\ast(\beta))$, which is zero in view of the
homotopy invariance 1.1 (e).

\medskip\noindent We have proved that the
pairing factors through $Pic(\overline C,C_\infty) \otimes V(C)$,
where
$$
Pic(\overline C,C_\infty)=Div(C)/\{f\in L(C)^\times\mid f=1\hbox{ on }
C_\infty\}
$$
is the divisor class group of modulus $C_\infty$, where the closed subscheme
$C_\infty\subset\overline C$ is identified with the corresponding Cartier
divisor and hence with the Weil divisor $\mathop{\sum}\limits_{y\in C_\infty}y$.

Let $Pic^0(\overline C,C_\infty)$ be the subgroup defined by the divisors
of degree zero, where $deg \ x = [\kappa(x):L]$. Then we have a canonical isomorphism
$$
Pic^0(\overline C,C_\infty)=J_{C_\infty}(\overline C)(L) \ ,
$$
where the right hand side is the group of $L$-rational points of the Rosenlicht generalized Jacobian
$J_{C_\infty}(\overline C)$ (cf. \cite{Ro}, \cite{Se} Chap. V and \cite{CC}),
which is a commutative, smooth, geometrically connected group variety
over $L$. Now
$$
J_{C_\infty}(\overline C)\times_L K = J_{C_{\infty,K}}(\overline C_K)=
J_{\tilde C_\infty}(\tilde C_K) \ , \leqno(1.2.1)
$$
where $C_{\infty,K} = C_\infty\times_L K =
C_\infty\times_{\overline C} \overline C_K\subseteq \overline C_K
= \overline C\times_L K$ as a closed subscheme, $\tilde C_K$ is
the complete regular model of the smooth curve $C_K=C\times_L K$,
and $\tilde C_\infty = C_\infty\times_{\tilde C}\tilde C_K$ as a
closed subscheme (Note that $\overline C_K$ is not necessarily
regular and that $C_{\infty,K}$ is not necessarily reduced, if
$char \ K>0$ and $K/L$ is not separable). In the language of
schemes the equalities (1.2.1) can most easily be seen by the fact that $J_D(X)$,
for a geometrically integral curve $X$ over a field $k$ and a
Cartier divisor $D\subseteq X$ such that $X-D$ is smooth,
represents the Picard functor $Pic_{X_D/k}$, where $X_D$ is the
curve obtained by contracting $D$ to a point, i.e., $X_D$ is the
scheme theoretic amalgamated sum $X\amalg_D Spec \ k$ (cf. [Se] p.
85 and [CC]). In fact, one has $(\overline C_K)_{C_{\infty,K}} =
(\tilde C_K)_{\tilde C_\infty}$, since the diagram
$$
\xymatrix{ \tilde C_\infty\ar[r]\ar[d] & \tilde C_K\ar[d]\\
C_{\infty,K}\ar@{^{(}->}[r]  & \overline C_K}
$$
is cartesian and cocartesian.

\medskip If the pairing
$$
\psi_K:Div(C_K)\times V(C_K)\longrightarrow V(K)
$$
is defined for $C_K$ over $K$ in the same way as for $C$ over $L$ above, then
by the same argument as there, this pairing factors through
$Pic(\tilde C_K,(\tilde C_\infty)_{red})\otimes V(C_K)$ and therefore also
through $Pic(\tilde C_K,\tilde C_\infty)\otimes V(C_K)$, where
$$
Pic(\tilde C_K,\tilde C_\infty)=Div(\tilde C_K)/\{f\in K(C_K)^\times\mid
f\equiv 1 mod \ \tilde C_\infty\}
$$
is the divisor class group of modulus $\tilde C_\infty~~ (f\equiv
1 mod \ \tilde C_\infty$ meaning that $f$ lies in
$\Gamma(U,\mathcal O_U)$ for an open neighbourhood of $(\tilde
C_\infty)_{red}$ and has image 1 in $\Gamma(\tilde C_\infty,
\mathcal O_{\tilde C_\infty})$). Let $p: C_K\to C$ be the
projection induced by $p:Spec \ K\to Spec \ L$. Then by 1.1 (b) we
have a commutative diagram of pairings
$$
\begin{array}{rccccc}
\psi_K: & Div(C_K) & \times & V(C_K) & \longrightarrow & V(K) \\
& \scriptstyle p^\ast\uparrow & & \uparrow \scriptstyle p^\ast & &
\uparrow \scriptstyle p^\ast
\\ \psi : & Div(C) & \times & V(C) & \longrightarrow & V(L) \, ,
\end{array}
$$
where the left hand map is the pull-back of Cartier divisors
(which sends $x\in\mid C\mid$ to the unique $x'\in\mid C_K\mid$
with $p(x')=x$; note that $\kappa(x)\otimes_L K$ is a field, since
$\kappa(x)$ is separable over $L$ and $L$ is separably closed in
$K$).

\medskip\noindent For the following we may assume that $V$ is
annihilated by a natural number $n$, since it
suffices to prove Theorem 1.2 for all subfunctors $V[n]$ for
such $n$. Then we get an induced diagram
$$
\begin{array}{cccccc}
 Pic(\tilde C_K,\tilde C_\infty)/n & \otimes & V(C_K) & \longrightarrow & V(K) \\
 \scriptstyle p^\ast\uparrow & & \uparrow \scriptstyle p^\ast & &
\uparrow \scriptstyle p^\ast
\\  Pic(\overline C,C_\infty)/n & \otimes & V(C) & \longrightarrow & V(L) \,
\end{array} ,\leqno(1.2.2)
$$
To prove Claim 1 in this situation it then suffices to show

\medskip\noindent {\bf Claim} 3: \ $p^\ast: \ Pic^0(\overline
C,C_\infty)/n\longrightarrow Pic^0 (\tilde C_K,\tilde
C_\infty)/n$ is surjective.

\medskip\noindent In fact, consider an element $\beta\in V(C)$ mapping
to our element $\alpha\in Im(V(C)\to V(K))$. Then $\alpha =
\psi_K(y_0,p^\ast(\beta))$, where $y_0$ is the $K$-rational point
of $C_K$ corresponding to the generic point $Spec \ K\to Spec \
F\to C$. Since $K/K^0$ is an $n$-rigid extension,
and $L$ is algebraically closed in $K$ by definition, there exists
an $L$-rational point $x_0$ of $C$.

\medskip\noindent
Then $y_0-p^\ast(x_0)$ lies in $Div^0(C_K)$, and
by Claim 3 there exists an element $z\in Div(C)$ with
$\psi_K(y_0-p^\ast(x_0),p^\ast(\beta))=\psi_K
(p^\ast(z),p^\ast(\beta))$. Thus $\alpha =
\psi_K(y_0,p^\ast(\beta)) = p^\ast\psi(x_0+z,\beta)$ lies in the
image of $p^\ast:V(L)\to V(K)$ as wanted.

\medskip\noindent To prove Claim 3, note that the considered map can be
identified with the natural map
$$
J_{C_\infty}(C)(L)/n\longrightarrow J_{C_\infty}(C)(K)/n \ , \leqno(1.2.3)
$$
since there are identifications $J_{C_\infty}(C)(K)=(J_{C_\infty}(C)\times_L K)(K)=
J_{\tilde C_\infty}(\tilde C_K)(K)= Pic(\tilde C_K,\tilde C_\infty)$, functorially in
the field $K$, and since $J_D(X)$ represents $Pic_{X_D/k}$ as mentioned above.
Since $K/K^0$ is an $n$-rigid extension, the map (1.2.3) is an isomorphism
by definition.

\medskip\noi
\bigskip\noindent It remains to prove the injectivity of the map $V(K^0)\to
V(K)$. Since $K=\mathop{\lim}\limits_{\longrightarrow} A_i$ for
smooth $K^0$-algebras $A_i$, by 1.1 (f) it suffices to show
that
$$
q^\ast_i:V(Spec \ K^0)\longrightarrow V(X_i)
$$
is injective for every $q_i:X_i=Spec \ A_i\to K^0$. But every
$X_i$ has a $K$-rational point, hence a $K^0$-rational point
$s_i:Spec \ K^0\to X_i$ by Lemma 1.3. Since $s_i^\ast
q_i^\ast=(q_i s_i)^\ast=id$, $q_i^\ast$ must be injective.
This finishes the proof of Theorem 1.2.

\medskip\noi
\begin{prop}\label{Proposition 1.3}
The following field extensions are $n$-rigid (see Definition 0.2).

\smallskip\noi (i) An extension $K/K^0$ of algebraically closed fields, where $n$ is arbitrary.

\smallskip\noi(ii) $K = \r$ and $K^0$ is dense and algebraically closed in $\r$, and $n$ is arbitrary.

\smallskip\noi(iii) $K$ is a complete discrete valuation field, $K^0$ is algebraically closed in $L$
and dense in $L$ for the valuation topology, and $n$ is invertible in $L$.
\end{prop}

\bigskip\noindent  \textbf{Proof} Let $L$ with $K^0\subset L \subset K$
be an intermediate field which is algebraically closed in $K$, and let $C$ be smooth curve
over $L$. We have to prove the properties 0.2 (a) and (b).

(a): Assume that $C$ has a $K$-rational point.
(i): If $K$ is algebraically closed, and $L$ is algebraically closed in $K$, then $L$ is algebraically closed as well,
so that $C$ has an $L$-rational point.
(ii) and (iii): If $K = \r$, then $C$ also has an $L$ rational point by Lemma 1.4 below.

(b): Let $\overline{C}$ be a proper regular model of $C$, and let $C_\infty \subset \overline{C}$
be a not necessarily reduced closed subscheme.
(i): If $K$ is algebraically closed, then $L$ is algebraically closed as well, so that the closed subscheme $C_\infty$
defined in the proof of Theorem 1.3 is reduced. Then $J_{C_\infty}(\overline{C})$ is a semiabelian variety,
so that $J_{C_\infty}(K)$ is divisible and $J_{C_\infty}(K)/n = 0$.
(ii) and (iii): In these two cases
$nJ_{C_\infty}(C)(K)$ is open in $J_{C_\infty}(C)(K)$ for the strong topology on this set,
i.e., the topology coming from the topology of $K$), since the morphism
$n: \ J_{C_\infty}(C)\to J_{C_\infty}(C)$ is \'{e}tale, $n$ being
invertible in $L$. Therefore the claim follows from

\begin{lem}\label{Lemma 1.3} \it Let $K$ be $\mathbb R$,
$\mathbb C$, or a complete discrete valuation field, let $L$ be a
dense subfield, and let $X$ be a scheme of finite type over $L$.
If $K$ is separable over $L$ or if $X$ is smooth over $L$, then
$X(L)$ is dense in $X(K)$ for the strong topology.
\end{lem}

In fact, for $x\in J_{C_\infty}(C)(K)$, $x + nJ_{C_\infty}(C)(K)$ is open,
hence by Lemma 1.3 there is $y\in J_{C_\infty}(C)(L)$ contained in $x + nJ_{C_\infty}(C)(K)$.

\medskip\noindent Lemma 1.3 is proved in \cite{KaSa}, Lemma 4, where the first case is
explicitly stated and reduced to the second case in the proof, and
where the reader may also find a definition of the strong topology
(called the usual topology there).

\medskip\noi
then $L$ is algebraically closed
as well (since $K$ is algebraically closed in $K$).
Hence there exists an $L$-rational point $x_0$ of $C$.
The same holds for the other two cases by Lemma 1.3 below,
since $C$ has the $K$-rational point $y_0$.

\begin{coro}\label{Corollary 1.4} If $K,K^0$ and $V$ are
as in Theorem 1.2, then for every scheme $X$ of finite type
over $K^0$ the restriction map
$$
V(X)\longrightarrow V(X_K)
$$
is an isomorphism.
\end{coro}

\medskip\noindent Indeed, the functor $V_X$ with $V_X(Y)=V(X\times_{K^0}Y)$
is again a rigid functor on the category of all noetherian
$K^0$-schemes.

\begin{rem}\label{Remark 1.5} Let $R$ be an excellent discrete
valuation ring, let $\hat R$ be its completion, and denote by $k$
and $K$ the fraction fields of $R$ and $\hat R$, respectively.
Then by definition (\cite{EGA IV}(2),7.8.2), $K$ is separable over $k$.
Hence, if $K^0$ is the algebraic (= separable) closure of $k$
in $K$, then the assumptions of Theorem 1.1 are fulfilled for $K$
and $K^0$. Note that $K^0$ is the fraction field of the
Henselization $\tilde R$ of $R$ \cite{EGA IV}(4),18.9.3). In
particular, if $K$ is a global field and $v$ is a place of $K$,
then theorem 1.1 holds for the pair $(K_{(v)},K_v)$, where $K_v$
is the completion of $K$ at $v$ and $K_{(v)}$ is the algebraic
closure of $K$ in $K_v$. In fact, if $v$ is non-archimedean, then
the corresponding valuation ring is excellent (\cite{EGA IV}(2),7.8.3(ii),(iii).
If $K_s$ is a separable closure of $K$ and
$w$ is a place of $K_s$ extending $v$, then $K_{(v)}$ can also be
replaced by the isomorphic decomposition field $K_s^{G_w}$, where
$G_w=\{\sigma\in Gal(K_s/K) \mid\sigma w=w\}$ is the decomposition
group of $w$, which is the classical Henselian field associated to
$K$ and $w$.
\end{rem}

\section{Examples of rigid functors}

The first two examples of rigid functors are given by \'{e}tale cohomology and algebraic $K$-theory.

\begin{prop}\label{Proposition 2.1} Let $\mathcal S$ be a
category of quasi-compact schemes, and let $\mathcal F$ be an
\'{e}tale torsion sheaf on $\mathcal S$ whose torsion is
invertible in $\mathcal S$. Assume that $\mathcal F\mid_X =
f^\ast\mathcal F\mid_Y$ for every morphism $f:X\to Y$ in $\mathcal
S$. (e.g., for any natural number $n$ we may take $\mathcal S =
Sch^{qc}/\mathbb Z[1/n]$, the category of quasi-compact schemes on
which $n$ is invertible, and $\mathcal F = \mathbb
Z/n(j)=\mu_n^{\otimes j}$, the $j$-th Tate twist of the constant
sheaf $\mathbb Z/n$, cf. \cite{Mi},  p. 163). Then for any integer
$i\geq 0$ the functor
$$
V(Y)=H^i(Y,\mathcal F)
$$
given by the $i$-th \'{e}tale cohomology with coefficients in
$\mathcal F$ is a rigid functor on $\mathcal S$.
\end{prop}

\medskip\noindent \textbf{Proof} \ This is well-known: The
contravariance of $V$ is induced by the adjunction maps
$ad_f:\mathcal F\to f_\ast f^\ast\mathcal F$ for morphisms $f:X\to
Y$ (cf. \cite{Mi} III 1.6(c)), viz., $f^\ast$ is the composition
$$
f^\ast:H^i(Y,\mathcal F) \ \mathop{\longrightarrow}\limits^{ad_f}
H^i(Y,f_\ast f^\ast\mathcal F) \
\mathop{\longrightarrow}\limits^{can}H^i(X,f^\ast\mathcal
F)=H^i(X,\mathcal F) \ .
$$
For the limit property 1.1(f) cf. \cite{Mi} III 1.16, and for the
homotopy property 1.1 (e), which is related to smooth base change,
cf. \cite{Mi} VI 4.15. The transfer $\pi_\ast$ for a finite flat
morphism $\pi:X\to Y$ is defined as follows. By \cite{SGA 4.3} XVII 6.2.3
there is a canonical trace morphism
$$
Tr_\pi \ : \ \pi_\ast\pi^\ast\mathcal F \ \longrightarrow \
\mathcal F
$$
for every abelian sheaf $\mathcal F$ on $Y$ (note that $\pi_\ast =
\pi_!$ for a finite morphism). Then $\pi_\ast$ is defined as the
composition
$$
\pi_\ast:H^i(X,\pi^\ast\mathcal F) \
{\mathop{\longleftarrow}\limits^{can}_\sim} \
H^i(Y,\pi_\ast\pi^\ast\mathcal F) \
\mathop{\longrightarrow}\limits^{Tr_\pi} \ H^i(Y,\mathcal F)
$$
in which the first map is an isomorphism since $\pi$ is finite (so
that $\mathcal G\mapsto\pi_\ast\mathcal G$ is exact, cf. \cite{Mi} II
3.6). The functoriality $(\rho\pi)_\ast = \rho_\ast\pi_\ast$ in
1.1 (a) then follows from the transitivity of the trace morphism
(\cite{SGA 4.3} XVII 6.2.3 ($Var 3$)). For property 1.1 (d) we first note
that the composition
$$
\mathcal F \ \mathop{\longrightarrow}\limits^{ad_\pi} \
\pi_\ast\pi^\ast\mathcal F \
\mathop{\longrightarrow}\limits^{Tr_\pi}\mathcal F
$$
is multiplication by $d$ for {\it every} finite flat morphism
$\pi: X\to Y$ of (constant) rank $d$ (loc.cit. (Var 4)). Hence
$\pi_\ast \pi^\ast$ is multiplication by $d$ for such $\pi$. Now
let $\pi:X_d= X[T]/(T^{d+1})\to X$ be as in 1.1 (d), and let
$s:X\to X_d$ be the obvious ``zero'' section $(T \mapsto 0)$ of
$\pi$. Then $s^\ast $ is an isomorphism by the topological
invariance of \'{e}tale cohomology (\cite{SGA 4.3} VIII 1.2), hence
$\pi^\ast$ is an isomorphism as well. Thus the equality
$\pi^\ast\pi_\ast \pi^\ast=d \pi^\ast$ implies that $\pi^\ast
\pi_\ast=d$ as wanted. If
$$
\xymatrix{ X' \ar[r]^{f'}\ar[d]_{\pi '} & X\ar[d]^{\pi}\\
Y'\ar[r]^{f} & Y }
$$
is a cartesian
diagram, with $\pi$ (and hence $\pi'$) finite and flat, then the
base change property $f^\ast\pi_\ast=\pi'_\ast f'^\ast$ of 1.1 (b)
follows from the commutativity of the diagram
$$
\xymatrix{ \pi_\ast f'_\ast f^{'\ast}\pi^\ast\mathcal
G\ar@{=}[d]^{\hspace{1,5cm}\displaystyle(1)} & &
\pi_\ast\pi^\ast\mathcal G\ar[rr]^{Tr_\pi}\ar[ll]_{\pi_\ast
ad_{f'}}\ar[d]^{ad_f \hspace{1cm}\displaystyle(2)} & & \mathcal
G\ar[d]^{ad_f}\\f_\ast \pi'_\ast f^{'\ast}\pi^\ast\mathcal G
\ar@{=}[d] & & f_\ast f^\ast \pi_\ast \pi^\ast \mathcal
G\ar[ll]_{f_\ast\Phi}\ar[rr]^{f_\ast f^\ast
Tr_\pi}\ar@{}[d]^{\hspace{-4mm}\displaystyle(3)} & & f_\ast
f^\ast\mathcal G\ar@{=}[d]\\f_\ast \pi'_\ast\pi^{'\ast}
f^\ast\mathcal G\ar[rrrr]^{f_\ast Tr_{\pi'}} & &  & & f_\ast
f^\ast\mathcal G}
$$
for any torsion sheaf $\mathcal G$ on $Y$, where $\Phi:
f^\ast\pi_\ast\to\pi'_\ast f'^\ast$ is the base change morphism.
Here (3) is commutative by loc. cit. (Var.2), (2) commutes by the
functoriality in sheaves of $ad_f$, and (1) commutes by the
definition of the base change morphism (cf., e.g, \cite{Mi} p. 223).
Note that we get an induced commutative diagram
$$
\xymatrix{ H^{i}(Y,\pi_{\ast}\pi^{\ast}\mathcal
G)\ar[r]\ar[d]^{\wr}_{can} & H^{i}(Y,\pi_{\ast} f'_{\ast}
f'^{\ast}\pi^{\ast}\mathcal G)\ar[d]^{\wr}_{can} & = &
H^{i}(Y,f^{\ast}\pi'_{\ast}\pi'^{\ast} f^{\ast}\mathcal
G)\ar[d]^{\wr}_{can}\ar[r] & H^{i}(Y,f_{\ast} f^{\ast}\mathcal
G)\ar[d]^{\wr}_{can}\\
H^{i}(X,\pi^{\ast}\mathcal
G)\ar[r]^{\hspace{-0,5cm}ad_{f'}}\ar[dr]_{f'^{\ast}} &
H^{i}(X,f'_{\ast} f'^{\ast}\pi^{\ast}\mathcal G)\ar[d]_{can} &&
H^{i}(Y',\pi'_{\ast}\pi'_{\ast}\pi'^{\ast} f^{\ast}\mathcal
G)\ar[r]^{\hspace{0,7cm}Tr_{\pi'}}\ar[d]^{\wr}_{can} &
H^{i}(Y',f^{\ast}\mathcal G)\\ &
H^{i}(X',f^{'\ast}\pi^{\ast}\mathcal G) & = & H^{i}(X',\pi^{'\ast}
f^{\ast}\mathcal G)\ar[ur]_{\pi '_{\ast}} & }
$$
Finally, property 1.1 (c) is a straightforward consequence of
\cite{SGA 4.3} XVII 6.2.3.1: it implies that for $X=X_1\amalg X_2$ with open and
closed immersions $\alpha_i:X_i\hookrightarrow X$ $(i=1,2)$ one
has $(\alpha_i)^\ast(\alpha_j)_\ast = \delta_{ij} id_{X_i}$ .

\begin{prop}\label{Proposition 2.2}\ For every integer $m\geq
0$, the functor
$$
V(X)=K_m(X)
$$
given by the $m$-th algebraic $K$-group is a rigid functor on the
category $\mathcal S$ of all noetherian schemes, except that the
homotopy axiom 1.1 (e) possibly only holds for a regular
base $X$.
\end{prop}

\medskip\noindent \textbf{Proof} \ This follows from the results of
Quillen in \cite{Qui}: By \cite{Qui} \S7,2, $X\mapsto K_m(X)$ is a
contravariant functor, and the limit property 1.1(f) is proved in
loc. cit. \S7, 2.2. The transfer map for a finite flat morphism
$\pi: X\to Y$ is defined as follows. Recall that $K_m(X) =
K_m(P(X))$, the $m$-th $K$-group of the exact category $P(X)$ of
locally free coherent $\mathcal O_X$-modules. Then $\pi_\ast:
K_m(X)\to K_m(Y)$ is induced by the exact functor $\pi_\ast:
P(X)\to P(Y)$ sending an $\mathcal O_X$-module $\mathcal P$ in
$P(X)$ to the $\mathcal O_Y$-module $\pi_\ast \mathcal P$ in
$P(Y)$ (cf. loc. cit. \S7). The functoriality $(\rho\pi)_\ast =
\rho_\ast\pi_\ast$ is immediate. For 1.1 (b) recall that the
pull-back $f^\ast: K_m(Y)\to K_m(Y')$ for a morphism $f:Y'\to Y$
is induced by the exact functor $f^\ast: P(Y)\to P(Y')$ sending
$\mathcal Q$ in $P(Y)$ to $f^\ast \mathcal Q = \mathcal
O_{Y'}\otimes_{f^{-1}\mathcal O_Y} f^{-1} \mathcal Q$ (coherent
pull-back) in $P(Y)$. If now
$$
\xymatrix{ X' \ar[r]^{f'}\ar[d]_{\pi '} & X\ar[d]^{\pi}\\
Y'\ar[r]^{f} & Y }
$$
is a cartesian diagram, with $\pi$ finite and flat, then the base
change morphism $f^\ast\pi_\ast\to\pi'_\ast f'^\ast$ is an
isomorphism of exact functors from $P(X)$ to $P(Y')$; therefore
$f^\ast\pi_\ast = \pi'_\ast f'^\ast$ on the level of $K$-groups
(\cite{Qui} \S1 Prop. 2). The additivity property 1.1 (c) follows
from \cite{Qui} \S1, (8) and the fact that for $X=X_1\amalg X_2$
the immersions $\alpha_i:X_i\hookrightarrow X$ induce an
equivalence of exact categories
$(\alpha^\ast_1,\alpha_2^\ast):P(X)\to P(X_1)\times P(X_2)$ with
$(\alpha_i)^\ast(\alpha_j)_\ast\stackrel{\sim}
{\longrightarrow}\delta_{ij}id \ : \ P(X_j)\to P(X_i)$. Property
1.1 (d) follows from the more general fact that for any flat
finite morphism $\pi: X\to Y$ for which $\pi_\ast\mathcal O_X\cong
\mathcal O^d_Y$ as an $\mathcal O_Y$-module one has a functorial
isomorphism
$$
\pi_\ast\pi^\ast P \stackrel{\sim}{\longrightarrow}
P\otimes_{\mathcal O_Y}\pi_\ast \mathcal O_X
\stackrel{\sim}{\longrightarrow} P^d
$$
for $P$ in $P(Y)$. This implies that $\pi_\ast\pi^\ast$ is the
multiplication by $d$ on $K_m(Y)$ by \cite{Qui} \S3 Corollary 1.
Finally, the homotopy property 1.1(e) for $V(X)=K_m(X)$ is proved
in \cite{Qui} \S7, Proposition 4.1, for a regular base $X$.

\smallskip\noindent
This covers case (2) in the introduction, while the following implies case (3) in the introduction.

\begin{theo}\label{Theorem2.3} The functors $X \mapsto CH^m(X,n)$ of Bloch's higher Chow groups are
rigid functors. In particular, this holds for the classical Chow groups $CH^r(X) = CH^r(X,0)$.
\end{theo}

\v\noi\textbf{Proof} All properties in Definition 1.1 hold, see for example \cite{Lev}, section 2.1.

\section{Sufficiently rigid functors}

\v\noi
In view of the applications we have in mind,
we note that the full strength of the axioms in 1.1 was not needed
in the proof of Theorem 1.2. Consider the following weakening.

\begin{defi}\label{Definition 3.1} Let $\mathcal S=Sch^{noeth}/k$ be the
category of noetherian $k$-schemes, for a field $k$. A
contravariant functor $V:\mathcal S\mapsto Ab$ is called sufficiently
rigid, if it satisfies the following properties.

\v\noi (a) Call a morphism $\pi:X\rightarrow Y$ in $\mathcal
S$ admissible, if it is a finite and flat morphism of smooth
$L$-schemes, for a field extension $L$ of $k$. Then for any
admissible morphism $\pi:X\rightarrow Y$ there is a transfer map
$\pi_\ast: V(X)\rightarrow V(Y),$ such that
$(\varrho\pi)_\ast=\varrho_\ast \pi_\ast$ for another admissible
morphism $\varrho:Y\rightarrow Z.$

\v\noi (b) For every cartesian diagram of schemes
$$
\xymatrix {X'\ar[r]^{f'}\ar[d]_{\pi '} & X\ar[d]^{\pi} \\
Y'\ar[r]^{f} & Y}
$$
lying in $\mathcal S$, with $\pi$ and $\pi'$ admissible, one has
$f^\ast\pi_\ast= \pi'_\ast f'^\ast:V(X)\to V(Y').$

\v\noi(c) Let $\pi:C\to D$ be a finite flat morphism of smooth
curves over an extension $L$ of $k$. For $y\in D$ and $x\in C$
with $\pi(x)=y$ consider the commutative (in general not
cartesian!) diagram
$$
\xymatrix{ Spec~(\kappa(x))\ar[r]^{\hspace{0,8cm}\varphi_x}\ar[d]_{\pi_x} & C\ar[d]^{\pi}\\
Spec~(\kappa (y))\ar[r]^{\hspace{0,8cm}\varphi_y} & D} \, \leqno(3.1.1)
$$
where $\varphi_x,\varphi_y$ and $\pi_x$ are the canonical
morphisms, and for every smooth $L$-scheme $X$ denote by
$$
\xymatrix{ X \times_L\kappa(x)\ar[r]^{\varphi_x}\ar[d]_{\pi_x} & X \times_L C\ar[d]^{\pi}\\
X \times_L\kappa(y)\ar[r]^{\varphi_y} & X \times_L D}
\leqno(3.1.2)
$$
also the diagram obtained by base change with $X$ over $L$. Then
$$
\varphi^\ast_y \pi_\ast=\sum_{\pi(x)=y}
e(x/y)(\pi_x)_\ast~(\varphi_x)^\ast :V(X \times_L C)\rightarrow V(X
\times_L\kappa(y)),
$$
where $e(x/y)$ is the ramification index of $x$ over $y$, and
the same equality holds after base change with an open subscheme
$U$ of $X \times_L D.$

\v\noi(d) If $X$ is a smooth $L$-scheme, for an extension $L$ of
$k$, then the projection $p:\mathbb A_X^1\rightarrow X$ induces an
isomorphism $p^\ast:V(X) \stackrel{\sim}{\rightarrow} V(\mathbb
A_X^1).$

\v\noi(e) Let $i\mapsto X_i$ be a filtered projective system of
schemes in $ \mathcal S$ such that the transition morphisms
$X_i\rightarrow X_j$ are affine, and assume that
$X=\mathop{\lim}\limits_{\longleftarrow} X_i$ is in $\mathcal S$.
Then the canonical map
$$
\lim_{\longrightarrow} V(X_i)\rightarrow V(X)
$$
is an isomorphism.
\end{defi}

Then we have

\begin{theo}\label{Theorem 3.2} Let $K/K^0$ be an $n$-rigid field extension (see Definition 0.2),
and let $V:Sch^{noeth}/K^0\rightarrow Ab$ be a
sufficiently rigid functor. Then for every smooth $K^0$-scheme
$X$ the restriction map
$$
V(X)\rightarrow V(X_K)
$$
is an isomorphism, if $V$ has values in abelian groups annihilated by $n$.
\end{theo}

\medskip\noindent\textbf{Proof} The proof of Theorem 1.2 applies to
the functor $V_X$ with $V_X(Y)=V(X \times_{K^0} Y).$ In fact,
by directly using 1.6 (c), applied to $X \times_{K^0}
C\stackrel{id_X\times\pi}{\longrightarrow} X
\times_{K^0}\mathbb A^1_L$, we need only consider transfer
maps for admissible morphisms, and we only need the homotopy axiom
1.1 (e) for $\mathbb A^1_{X
\times_{K^0}L}=X\times_{K^0}\mathbb A^1_L\rightarrow
X\times_{K^0}L$, with $L$ and $\pi:C\rightarrow\mathbb A^1_L$
as in the proof of theorem 1.2. Note that $X\times_{K^0}
Y=(X\times_{K^0}L)\times_L Y$ is a smooth $L$-scheme for a
smooth $L$-scheme $Y,L$ an extension of $K^0.$

\medskip Of course, here we only needed 3.1 (c) for the diagrams (3.1.1),
and not for their localizations. This additional property will be
needed in section 2. Therefore we record the following, not
completely obvious fact.

\begin{lem}\label{Lemma 3.3} Let $k$ be a field. A rigid
functor $V:Sch^{noeth}/k \rightarrow Ab$ is also a sufficiently
rigid functor. \

\medskip\noindent \textbf{Proof} We only have to show that 3.1 (c)
holds for $V$. With the notations of 3.1 (c) we have a commutative
diagram
$$
\xymatrix{ \scriptstyle\coprod\limits_{\pi(x)=y} \textstyle
X\times_L\kappa(x)
\ar[r]^{\hspace{0,5cm}\amalg\alpha_x}\ar[dr]_{\amalg\pi_x} &
X\times_L~C_y\ar[d]^{\pi'} \ar[r]^{\varphi'_y} & X\times_L
C\ar[d]^{\pi} \\ & X\times_L\kappa(y)\ar[r]^{\varphi_y} &
X\times_L D}
$$
such that the square is cartesian and
$\varphi'_y\alpha_x=\varphi_x.$ By assumption, the triangle can be
rewritten and factored as
$$
\xymatrix{\scriptstyle\coprod\limits_{\pi(x)=y}\textstyle
X\times_L\kappa(x) \ar[r]^{\hspace{1,3cm}(\alpha_x)} \ar[drr]_{id}
& \scriptstyle\coprod\limits_{\pi(x)=y}\hspace{-1cm}
& X\times_L\left(\kappa(x)[T]/(T^{e(x/y)})\right)\ar[d]^{(\gamma_x)} \\
& & \scriptstyle\coprod\limits_{\pi(x)=y}\textstyle
X\times_L\kappa(x)\ar[d]^{\amalg_{\pi_x}}
\\ &  & X\times_L\kappa(y).}
$$
By forming the base change of the diagrams with an open subscheme
$U$ of $X\times_L D,$ we obtain a commutative diagram
$$
\xymatrix{\scriptstyle\coprod\limits_{\pi(x)=y}\textstyle Y_x
\ar[r]^{\hspace{0,8cm}(\alpha_x)} \ar[drr]_{id} &
\scriptstyle\coprod\limits_{\pi(x)=y} \hspace{-1cm}&
Y_x[T]/(T^{e(x/y)+1}) \ar[r]^{\varphi'_y} \ar[d]^{(\gamma_x)} &
U'=U\times_D C
\ar[dd]^{\pi}\\
&  & \scriptstyle\coprod\limits_{\pi(x)=y}
Y_x\ar[d]^{\scriptstyle\coprod_{\pi_x}}\\
& & Y\ar[r]^{\varphi_y} & U & ,}
$$
with cartesian square, and we want to show that
$$
\varphi_y^\ast
\pi_\ast=\sum_{\pi(\times)=y}e(x/y)(\pi_x)_\ast\varphi^\ast_x:
V(U')\rightarrow V(Y),
$$
where $\varphi_x=\varphi'_y \alpha_x~.$ Since ${\varphi'}_y^\ast
\pi_\ast= (\amalg \pi_x)_\ast (\gamma_x)_\ast(\varphi_y')^\ast$ by
1.1 (b), it suffices to show that
$$
(\gamma_x)_\ast=e(x/y)\alpha^\ast_x,
$$
by 1.1 (c). Since $\gamma_x\alpha_x=id_{Y_x}$, we have
$\alpha^\ast_x \gamma^\ast_x=id$, and $\gamma^\ast_x$ is
injective. Hence it suffices to show the equality
$$
\gamma^\ast_x(\gamma_x)_\ast=e(x/y)\gamma^\ast_x \alpha^\ast_x=e(x/y)id,
$$
which holds by 1.1 (d).
\end{lem}

\section{The treatment of associated Zariski sheaves}

\v\noi
The following considerations will be useful for treating $\mathcal K$- and $\mathcal H$-cohomologies.

\v Let $\mathcal S $ be a category of schemes such that for a morphism
$f:X\to Y$ in $\mathcal S$ and an open immersion
$j:U\hookrightarrow Y$ all morphisms of the cartesian diagram
$$
\xymatrix{ U'\ar@^{^{(}->}[r]^{j'}\ar[d]_{f'} & X\ar[d]^{f}\\
U\ar@{^{(}->}[r]^{j} & Y}\leqno(4.1.1)
$$
lie in $\mathcal S$. Thus $\mathcal S$ becomes a site, if we endow
it with the Zariski topology. Let $V$ be a contravariant functor
on $\mathcal S$ with values in the category $Ab$ of abelian
groups, i.e., a presheaf (of abelian groups) on $\mathcal S$, and
let $\mathcal V$ be the associated Zariski sheaf on $\mathcal S.$
If $V_X$ and $\mathcal V_X$ denote the restrictions of $V$ and
$\mathcal V$ to the small Zariski site $X_{Zar}$ for a scheme $X$
in $\mathcal S$ ( which consists of all open immersions
$(U\hookrightarrow X)$, then $\mathcal V_X= a~V_X,$ where $a=a_X$
is the functor mapping a presheaf on $X_{Zar}$ to its associated
sheaf (on $X_{Zar}$).

\medskip For a scheme $X$ in $\mathcal S$,
let $H^i(X,\mathcal V)$ be the i-th Zariski cohomology of
$\mathcal V$ on $X$. This is equivalently, the i-th derived
functor of the functor $\mathcal F\mapsto H^0(X,\mathcal
F)=\mathcal F (X),$ from sheaves on $\mathcal S$ to $Ab$,
evaluated at $\mathcal V$, and also equal to $H^i(X,\mathcal
V_X),$ the i-th cohomology of $\mathcal V_X$ on $X_{Zar}$ (cf.
\cite{Mi} III 1.5 (b), 1.10, and 3.1 (c)). It is clear from the first
description that $X\shortmid\hspace{-0,8mm}\rightsquigarrow H^i
(X,\mathcal V)$ is a contravariant functor from $\mathcal S$ to
$Ab$. In the second description, this functoriality can be
described as follows. Since $\mathcal V$ is a sheaf on $\mathcal
S$, one gets a canonical morphism
$$
\alpha_f:\mathcal V_Y\to f_\ast \mathcal V_\ast
$$
for every morphism $f:X\to Y$ in $\mathcal S$, defined via the
maps
$$
\mathcal V_Y(U)=\mathcal
V(U)\stackrel{(f')^\ast}{\longrightarrow}\mathcal V(U')=\mathcal
V_X(U')=f_\ast \mathcal V_\ast(U)
$$
for each diagram (4.1.1). Then the wanted pull-back is the
composition
$$
H^i(Y,\mathcal
V_Y)\stackrel{\alpha_f}{\longrightarrow}H^i(Y,f_\ast\mathcal
V_\ast) \stackrel{can}{\longrightarrow} H^i(X,\mathcal
V_X).\leqno(4.1.2)
$$
\noi We have the following result.

\begin{theo}\label{Theorem 4.1} Let $\mathcal S=
Sch^{noeth}/k$ be the category of noetherian $k$-schemes, for a
field $k$, and let $V$ be a sufficiently rigid functor on
$\mathcal S$. Assume that the following properties hold for rings
$A$ which are localizations of (the coordinate rings of) smooth
affine $L$-schemes $Y$, where $L$ is an extension field of $k$.

\v\noi (i) If $A$ is semi-local, then the natural map $V(Spec ~A)\rightarrow H^0(Spec~A,\mathcal V)$
is an isomorphism, and $H^i(Spec~A,\mathcal V)=0$ for $i>0$.

\v\noi(ii) If $A$ is local, then the restriction map $H^0(Spec~A,\mathcal V)\to H^0(Spec A[T],\mathcal V)$
is an isomorphism, and $H^i(Spec~A[T], \mathcal V)=0$ for $i>0.$

\v\noi Then for every $\nu\geq 0$ the functor
$$
Y\quad\shortmid\hspace{-0,07cm}\rightsquigarrow\quad
H^\nu(Y,\mathcal V)
$$
is a sufficiently rigid functor on $\mathcal S.$
\end{theo}

\medskip\noindent \textbf{Proof} For a scheme $Y$ denote by
$P(X_{Zar})$ and $S(X_{Zar})$ the category of abelian pre\-sheaves
and sheaves on the small Zariski site $X_{Zar}$, respectively, let
$i: S(X_{Zar})\to P(X_{Zar})$ be the inclusion, and let $a=a_X:
P(X_{Zar} )\to S(X_{Zar})$ be its left adjoint, mapping a presheaf
$P$ to its associated sheaf. For a morphism $f:X\rightarrow Y$ let
$f_P:P(X_{Zar})\to P(Y_{Zar})$ be the direct image functor
(defined by $(f_P P)(U)=P(f^{-1}(U)))$, and let $f^P:
P(Y_{Zar})\to P(X_{Zar})$ be its left adjoint.

\medskip We first show the limit property 3.1 (e) for the functors
$H^\nu (-,\mathcal V)$. Let $(X_i)_{i\in I}$ be a filtered
projective system of schemes in $\mathcal S$ with affine
transition maps. Assume that $X=\mathop{\textstyle\lim}\limits_{\longleftarrow\, i} X_i$
is in $\mathcal S$, and let $u_i:X\to X_i$ be the canonical
morphism. Then the morphism
$$
\mathop{\textstyle\lim}\limits_{\longrightarrow\, i}~~ u_i^\ast\mathcal V_{X_i} \longrightarrow \mathcal V_X
$$
is an isomorphism. In fact, the stalk at a point $x\in X$ with
images $x_i$ in $X_i$ is the map
$$
\mathop{\textstyle\lim}\limits_{\longrightarrow\, i}~V(\mathcal O_{X_i, x_i})\to V(\mathcal O_{X,x})~,
$$
which is an isomorphism by 3.1 (e) for $V$, since $\mathcal
O_{X,x}= \mathop{\lim}\limits_{\longrightarrow\, i}  \mathcal O_{X_i,x_i}~.$

\medskip\noi It now follows that the map
$$
\lim_{\longrightarrow\,i} H^\nu(X_i,\mathcal V_{X_i})\to
H^\nu(X,\mathcal V_X)
$$
is an isomorphism since $X_{Zar}$ is the limit of the sites
$(X_i)_{Zar}$ (cf. \cite{SGA 4.2} VII 5.7 for the case of \'etale sites;
the case of Zariski sites is much simpler and follows from
\cite{EGA IV} 3, 8.6.3 and 8.10.5 (vi), cf. the proof of \cite{SGA 4.2} VII 5.6).

\medskip Now let $\pi:X\to Y$ be an admissible finite and flat
morphism in $\mathcal S$. The transfer morphisms for $\pi$ and all
base changes with open immersions $U\hookrightarrow Y$ define a
morphism
$$
\pi_P V_X\to V_Y~~.
$$
We get an induced diagram
$$
\xymatrix{ & & a~ \pi_P V_X \ar[r]\ar[d]^{\phi} &  a~V_Y  =  \mathcal V_Y\\
\pi_\ast\mathcal V_X & = & \pi_\ast a~V_X~,}
$$
where the morphism $\phi$ is the canonical map, which for instance
arises from the adjunction map $id\rightarrow i~a$ and the
equality $\pi_\ast = a~ \pi_P i$. We claim that $\phi$ is an
isomorphism, Let $y\in Y.$ Then the stalk of $\phi$ at $y$ is the
canonical map
$$
\xymatrix{ \mathop{\lim}\limits_{\longrightarrow\,U}~V(U\times_Y X)
\ar[r]\ar[d]_{\wr} &
\mathop{\lim}\limits_{\longrightarrow\,U} ~H^0(U\times_Y X,\mathcal V)\ar[d]^{\wr}\\
V({Spec}~ \mathcal O_{Y,y}\times_Y
X)\ar[r]^{\hspace{-0,5cm}\phi_y} & H^0(Spec~ \mathcal
O_{Y,y}\times_Y X, \mathcal V)~~,}
$$
where $U$ runs through all open neighborhoods of $y$ in $Y$. The
vertical maps are isomorphisms since 3.1 (e) holds for $V$ and
$H^0(-,\mathcal V).$ Since $\pi$ is finite, $Spec~\mathcal
O_{Y,y}\times_Y X=Spec~A$ for a semi-local ring $A$, which is a
localization of $X$. Since $\pi$ is admissible, $X$ is smooth over
an extension $L$ of $k$. Hence $\phi_y$ is an isomorphism by our
assumption (i) on $V$. As a consequence, we have a canonical
morphism
$$
Tr_\pi:~ \pi_\ast\mathcal V_X\rightarrow\mathcal V_Y~~.
$$
\smallskip\noi Next we claim that $R^i\pi_\ast~\mathcal V_X=0$ for $i>0$.
In fact, the stalk at  $y\in Y$ of $R^i\pi_\ast~\mathcal V_X$ is
$$
\lim_{\longrightarrow\,U}~H^i(U\times_Y X,\mathcal V),
$$
where the limit is over all open neighborhoods of $y$ in $Y$.
This equals $H^i(Spec~\mathcal O_{Y,y}\times_Y X,\mathcal
V)$ by the limit property 3.1 (e) for $H^i(-,\mathcal V)$, and
this is zero for $i>0$ by the assumption (i) on $V$ and the same
arguments as before.

\medskip As a consequence, the canonical map $H^\nu(Y,\pi_\ast
\mathcal V_X) \stackrel{can}{\rightarrow} H^\nu(X,\mathcal V_X)$
is an isomorphism for all $\nu\geq 0$ , and we define the transfer
maps for $\pi$ as the compositions
$$
\pi_\ast:H^\nu(X,\mathcal V_X)\stackrel{\sim}{\longleftarrow}
H^\nu(Y,\pi_\ast\mathcal V_X) \stackrel{Tr_\pi}{\longrightarrow}
H^\nu(Y,\mathcal V_Y)~~,
$$
for all $\nu\geq 0$. That
$(\varrho\pi)_\ast=\varrho_\ast\pi_\ast$, for an admissible
morphism $\varrho:Y\rightarrow Z$, is a straightforward
consequence of the fact that $Tr_{\varrho\pi}$ coincides with the
composition
$$
\varrho_\ast \pi_\ast \mathcal V_X \stackrel{\varrho_\ast
Tr_\pi}{\longrightarrow} \varrho_\ast \mathcal V_Y
\stackrel{Tr_\varrho}{\longrightarrow} \mathcal V_Z~~,
$$
which easily follows from the definitions and the observation that
$a \varrho_P \pi_P V_X\rightarrow \varrho_\ast a\pi_P V_X$ is an
isomorphism.

\medskip Now we show 3.1 (b). Let
$$
\xymatrix{ X'\ar[r]^{f'}\ar[d]_{\pi '} & X \ar[d]^{\pi}\\
Y'\ar[r]^f & Y}
$$
be a cartesian diagram in $\mathcal S$, with $\pi$ and $\pi'$
admissible. Then the diagram
$$
\xymatrix{ \pi_\ast~\mathcal V_X
\ar[r]^{\pi_\ast\alpha_{f'}}\ar[d]_{Tr_{\pi}} & \pi_\ast f'_\ast
\mathcal V_{X'} & = & f_\ast \pi'_\forall
\mathcal V_{X'}\ar[d]^{f_{\ast}Tr_{\pi'}}\\
\mathcal V_Y \ar[rrr]^{\alpha_f} & & & f_\ast \mathcal
V_{Y'}}\leqno(4.1.3)
$$
is commutative. In fact, by passing to the stalks at $y\in Y$ we
get the diagram
$$
\xymatrix{ V(Spec~ \mathcal O_{Y,y} \times_Y X)\ar[d]_{\pi_\ast}
\ar[r]^{f'^\ast}
& V(Spec~ \mathcal O_{Y,y}\times_Y X')\ar[d]^{\pi'_\ast}\\
V(Spec~ \mathcal O_{Y,y})\ar[r]^{f^\ast} & V(Spec~ \mathcal
O_{Y,y} \times_Y Y')}
$$
which is commutative by 3.1 (b) for $V$. The morphism $ \pi:
Spec~ \mathcal O_{Y,y}\times_Y X \rightarrow Spec~\mathcal
O_{Y,y}$ is only the localization of an admissible morphism, but
the map $\pi_\ast$ can be defined via passing to the limit over
the maps $\pi_\ast:V(U\times_Y X)\rightarrow V(U)$ for $U$ running
through the open neighborhoods of $y$ in $Y$. The analogous
statement holds for $\pi'_\ast$ and the equality
$f^\ast\pi_\ast=\pi'_\ast f'^\ast$ then follows from the
corresponding equalities for the $U$.

\medskip The equality $f^\ast\pi_\ast=\pi'_\ast f'^\ast$ for the
functors $H^\nu(-,\mathcal V)$ now follows from the commutative
diagram
$$
\xymatrix{ H^\nu(X,\mathcal V_X) \ar[r]^{\alpha_{f'}} &
H^\nu(X,f'_\ast\mathcal V_{X'}) \ar[rr]^{can}& & H^\nu(X',\mathcal
V_{X'})\\
H^\nu(Y,\pi_\ast\mathcal V_X)\ar[d]^{Tr_{\pi}}\ar[u]^{\wr}_{can}
\ar[r]^{\hspace{-5mm}\pi_\ast\alpha_{f'}}
 & H^\nu(Y,\pi_\ast f'_\ast \mathcal V_{X'})\ar@{}[d]^{\hspace{-5mm}\displaystyle(4.1.3)}
 \ar[u]_{can}\hspace{1cm}= & H^\nu(Y,f_\ast
\pi'_\ast \mathcal V_{X'})\ar[r]^{can}\ar[d]^{f_{\ast}Tr_{\pi}} &
H^\nu(Y',\pi'_\ast
\mathcal V_{X'})\ar[d]^{Tr_{\pi'}}\ar[u]^{\wr}_{can}\\
H^\nu(Y,\mathcal V_Y) \ar[rr]^{\alpha_f}  & & H^\nu(Y,f_\ast
\mathcal V_{Y'})\ar[r]^{can} & H^\nu(Y',\mathcal V_{Y'}).}
$$

\medskip The proof of 3.1 (c) for the $H^\nu(-,\mathcal V)$ is
similar. By similar arguments as above, it suffices to show that
the diagram
$$
\xymatrix{ \pi_\ast\mathcal V_{X\times C}\ar[d]_{Tr_{\pi}}
\ar[rr]^{\hspace{-2,2cm}\scriptscriptstyle\bigoplus \pi_\ast
\alpha_{\varphi_x}} & &
\stackrel{\scriptscriptstyle\bigoplus}{\scriptstyle\pi(x)=y}
\hspace{2mm}\pi_\ast(\varphi_x)_\ast~ \mathcal
V_{X\times\kappa(x)}\hspace{1cm}= &
\stackrel{\scriptscriptstyle\bigoplus}{\scriptstyle
\pi(x)=y}\hspace{1mm} (\varphi_y)_\ast)(\pi_x)_\ast~\mathcal
V_{X\times\kappa(x)}
\ar[d]^{\sum e(x/y)(\varphi_y)_\ast Tr_{\pi_\ast}}\\
\mathcal V_{X\times D} \ar[rrr]^{\alpha_{\varphi_y}} & & &
(\varphi_y)_\ast~\mathcal V_{X\times\kappa(y)}}
$$
commutes, where the notations are as in 3.1 (c). But by taking the
stalks at $t\in X\times \kappa(y)$ we obtain the diagram
$$
\xymatrix{ V(Spec~(R)\times_D C) \ar[d]_{\pi_\ast}
\ar[r]^{\hspace{-9mm}\scriptscriptstyle\bigoplus\varphi_x^\ast} &
\scriptstyle\bigoplus \limits_{\scriptstyle\pi(y)=x}
\displaystyle V\left( Spec~(R)\times_D \kappa(x)\right)\ar[d]^{\Sigma~e(x/y)(\pi_x)_\ast}\\
V(Spec~R) \ar[r]^{\hspace{-1cm}\varphi_y^\ast} & V\left(
Spec~(R)\times_D \kappa(x)\right)}
$$
with $R=\mathcal O_{X\times D, t}$, which is commutative by 1.6
(c) for $V$ (for $\pi_\ast$ and the $(\pi_x)_\ast$ the same remarks
as before apply).

\medskip Finally, we prove the homotopy property 3.1 (d) for the
functors $H^\nu(-,\mathcal V)$. Let $Y$ be smooth over an
extension $L$ of $k$, and let $p:\mathbb A_Y^1\rightarrow Y$ be
the affine line over $Y$. Then
$$
\alpha_p:\mathcal V_Y\rightarrow p_\ast\mathcal V_{\mathbb A^1_Y}
$$
is an isomorphism, and $R^i p_\ast \mathcal V_{\mathbb A_Y^1}=0$
for $i>0~,$ by assumption (ii) on $V$. In fact, for $y\in Y~,$ the
stalk of $\alpha_f$ at $y$ is the pull-back map
$$
H^0(Spec~\mathcal O_{Y,y},\mathcal V)\longrightarrow
H^0(\mathbb A^1_{Spec~\mathcal O_{Y,y}}, \mathcal V)~~,
$$
and the stalk $(R^i p_\ast\mathcal V)_y$ is isomorphic to $H^i
(\mathbb A^1_{Spec~\mathcal O_{Y,y}},\mathcal V)$ for all
$i\geqslant 0~.$ As a consequence, we obtain the bijectivity of
the maps in the composition
$$
p^\ast:H^\nu(Y,\mathcal
V_Y)~\raisebox{-6pt}{$\stackrel{\stackrel{\alpha_p}{\longrightarrow}}{\sim}$}~
H^\nu(Y,p_\ast \mathcal V_{\mathbb
A^1_Y})~\raisebox{-6pt}{$\stackrel{\stackrel{can}{\longrightarrow}}{\sim}$}~
H^\nu(\mathbb A_Y^1,\mathcal V_{\mathbb A_Y^1}),
$$
for every $\nu\geq 0~.$ q.e.d.

\begin{rem}\label{Remark 4.2} (a) It seems unlikely that one can define natural transfer
maps
$$
\pi_\ast:H^\nu(X,\mathcal V)\rightarrow H^\nu(Y,\mathcal V)
$$
for arbitrary finite and flat maps $\pi:X\rightarrow Y~,$ even if
one has transfer maps $\pi_\ast:V(X)\rightarrow V(Y)$ for all such
$\pi$.

\noi(b) Fix an extension $L$ of $k$, and let $p:\mathbb
A_Y^1\rightarrow Y$ be the affine line over a smooth $L$-scheme
$Y$. Given the limit property 1.6 (e) for $V$ (and hence for the
$H^\nu(-,\mathcal V))$ the following statements are equivalent.

\begin{itemize}

\item[(1)] $p^\ast: H^\nu(U,\mathcal V)\rightarrow H^\nu(\mathbb
A_U^1,\mathcal V)$ is an isomorphism for all $\nu\geqslant 0$ and
all open subschemes $U$ of $Y$.

\item[(2)] The morphism $\mathcal V_Y\rightarrow Rp_\ast\mathcal V_{\mathbb A^1_Y}$ is a quasi-isomorphism.

\item[(3)] Property 4.1 (ii) holds for all local rings $\mathcal
O_{Y,y}$ of $Y$.

\end{itemize}
\end{rem}

\section{$\cK$-cohomology and Gillet's Chow homology}

\begin{theo}\label{Theorem 5.1} Let $k$ be any field. The functor
$$
Y\mapsto H^\nu(Y,\mathcal K_m)~~,
$$
where $\cK_m$ is the Zariski sheaf associated to the presheaf
$U\mapsto K_m(U)$ given by the m-th algebraic $K$-groups, is a
sufficiently rigid functor on the category of all regular noetherian
$k$-schemes. If $F$ is a functor from abelian groups to abelian
groups, then the same holds for the functor $F\cK_m$ which is the
Zariski sheaf associated to the presheaf $U \mapsto F(K_m(U))$.
\end{theo}

\medskip\noindent \textbf{Proof} By Theorem 4.1 we have to show the properties 4.1
(i) and (ii). The first property is a direct consequence of
results of Quillen. In fact, let $Y$ be a noetherian scheme. In
[Qui] \S 7,5.4 Quillen constructed a spectral sequence
$$
E_1^{p,q}(Y)=\mathop{\oplus}\limits_{x\in Y^{(p)}}
K_{-p-q}(\kappa(x))\Rightarrow K'_{-p-q}(Y)~~,\leqno(5.1.1)
$$
which is contravariant for flat morphisms. Here $K'_m(Y)$ is the
m-th $K$-group of the category $M(Y)$ of coherent $\mathcal
O_Y$-modules, but for a regular noetherian scheme $Y$, this group
is canonically isomorphic to $K_m(Y)$ ([Qui] \S 7, (1.1)),
functorially for flat pull-backs (both functorialities are induced
by the exact functor mapping an $\mathcal O_Y$-module $\mathcal M$
to its coherent pull-back $f^\ast\mathcal M=\mathcal O_{Y'}\mathop{\otimes}_{f^{-1}\mathcal O_Y} f^{-1} \mathcal M$
for a flat morphism $f:Y'\rightarrow Y).$ If $Y$ is smooth over a
field $L$, or any localization of such a scheme, then Quillen
constructed a canonical isomorphism
$$
E_2^{p,q}(Y)=H^p(Y,\mathcal K_{-q})~~~(p,q\in\mathbb Z),\leqno(5.1.2)
$$
compatible with flat pull-backs, by showing that for any ring $A$
obtained by localizing (an open affine subscheme of) a smooth variety $Y$ in
finitely many points, the edge morphism
$$
K_{-m}(A)\longrightarrow E_2^{0,m}(Spec~ A)\leqno(5.1.3)
$$
is an isomorphism for all $m\in\mathbb Z$ and
$$
E_2^{p,q}(Spec~A)=0\quad\hbox{for }p\not= 0\leqno(5.1.4)
$$
(loc. cit. 5.6, 5.8 and 5.10). This proves 4.1 (i) for such $A$: by
definition, the composition
$$
K_{-m}(Spec~A)\rightarrow E_2^{0,m}(Spec~A)= H^0(Spec A, \mathcal K_{-m})\leqno(5.1.5)
$$
is the canonical map.

\medskip In view of (5.1.2) and Remark 4.2 (b), property 4.1 (ii)
follows from results of Gillet. For it follows from \cite{Gi} Thm.~8.3
that for the affine line $p:\mathbb A_Y^1\rightarrow Y$ over a
smooth $L$-scheme Y, the pull-back
$$
p^\ast:E_2^{p,q}(Y)\longrightarrow E_2^{p,q}(\mathbb A_Y^1)\leqno(5.1.6)
$$
is an isomorphism for all $p,q\in\mathbb Z$. In fact, without restriction,
$Y$ is equidimensional of dimension $d$, and then (4.3.6) is the map
$$
p^\ast: CH_{d-p,\,d+q}(Y)\rightarrow CH_{d+1-p,\,d+1+q}(\mathbb A_Y^1)\leqno(5.1.7)
$$
of Gillet's Chow homology groups in loc. cit., which is an isomorphism.

\medskip\noi
For the second claim of the Theorem we note that we obviously obtain sufficiently rigid
functors after applying the functor $F$. So we have to show the properties 4.1 (i) and 4.1 (ii).
For this we use a result of Grayson \cite{Gray} on the
universal exactness of the Gersten resolution for the $K$-theory of a smooth variety.
It implies that the statements remain true if one replaces $K_{-p-q}(k(x))$ by
$FK_{-p-q}(k(x))$ in (4.3.1), $\cK_{-q}$ by $F \cK_{-q}$ in (4.3.2), and $K_{-m}(A)$
by $FK_{-m}(A)$ in (4.3.3). Note that we do not assume that the functor $F$ is exact,
so that we replace $E_2^{p,q}(Y)$ by
$$
E_2^{p,q}(Y)_F := H^p(FE^{\ast,q}(Y)) \leqno(5.1.8)
$$
everywhere, which can be different from $FE_2^{p,q}(Y)$. But by Grayson's result (5.1.2) is then
replaced by an isomorphism
$$
E_2^{p,q}(Y)_F = H^p(Y,F \cK_{-q}) \,.\leqno(5.1.9)
$$
Again by universal exactness we get exactness of the complex
$$
0 \rightarrow FK_{-m}(A) \rightarrow FE^{0,m}_1(A) \rightarrow FE^{1,m}_1(A) \rightarrow \ldots
$$
for $A$ as above and hence the $F$-analog of (5.1.3), (5.1.4) and (5.1.5), and therefore property 4.1 (i).

\medskip\noi
As for the $F$-analog of 4.1 (ii), this is reduced, as in \cite{Gi} Theorem 8.3, to the case
of an affine line ${\mathbb A}^1_k \rightarrow \mbox{Spec}(k)$, and thus to an $F$-analog of \cite{Gi} Lemma 8.4.
But it is straightforward to check that all steps in the proof of that lemma go through for the
$F$-analog.

\begin{coro}\label{Corollary 5.2}
The claim of Theorem 0.3 is true for the case (3) of Gillet's Chow homology
groups.
\end{coro}

\medskip\noindent \textbf{Proof} Gillet's Chow homology groups of schemes are defined as
$$
CH_{r,s}(X) = E^2_{r,-s}(X)\leqno(5.2.1)
$$
for the $E^2$-terms of the spectral sequence
$$
E^r_{p,q}(X) \Rightarrow K'_{-p-q}(X)
$$
for the $K$-theory of the category $\cM(X)$ of quasi-coherent sheaves on $X$ and the filtration by dimension of support,
whose $E^1$-term is
$$
E^1_{p,q}(X) = \bigoplus_{x\in X_p} K_{p+q}(k(x)) \,,
$$
where $X_p$ is the set of points $x\in X$ of dimension $p$, i.e., whose closure is of this dimension,
and where $k(x)$ is the residue field of $x$.
If we consider varieties $X$ over a field $k$ and an extension field $K$ of $k$, then we have an exact functor from $\cM(X)$ to $\cM(X_K)$,
where $X_K=X\times_kK$. Since $Y$ and $Y_K$ have the same dimension for a $k$-variety $Y$, this induces a morphism of spectral sequences
$$
\xymatrix{
E^r_{p,q}(X_K)        & \Rightarrow &  K'_{-p-q}(X_K) \\
E^r_{p,q}(X) \ar[u]   & \Rightarrow &  K'_{-p-q}(X)  \ar[u] \,,
}
$$
and in particular morphisms
$$
E^2_{p,q}(X) \longrightarrow E^2_{p,q}(X_K)\,.
$$
The claim of Corollary 5.2 now follows from the following observation.

\begin{theo}\label{Theorem 5.3}
Let $\ell$ be a prime and assume that $K/k$ is an $\ell^m$-rigid field extension for all natural numbers $m$. Then the maps
$$
E^2_{p,q}(X)[\ell^n] \longrightarrow  E^2_{p,q}(X_K)[\ell^m] \hspace{1cm} \mbox{    and    } \hspace{1cm} E^2_{p,q}(X)/\ell^n \longrightarrow E^2_{p,q}(X_K)/\ell^m
$$
are isomorphisms for all $k$-varieties $X$, all natural numbers $m$ and all integers $p,q$, if $k$ is perfect or if $\ell$ is invertible in $k$.
\end{theo}

\medskip\noindent \textbf{Proof}
We proceed by induction on dimension.
First assume that $k$ is perfect. Then the claim is true for all $k$-varieties $X$ of dimension $0$, because these
are smooth and we have isomorphisms $E_2^{0,q}(X)=K_{-q}(X)=H^0(X,\cK_{-q})$, and the groups are zero for all
other degrees.

We use the last claim of Theorem 4.3 for the functor $\Tl$ on abelian groups, which send an abelian group $A$
to its $\ell$-primary torsion subgroup $\Tl(A)$, and for the functor $\Cl$, which sends $A$ to $\Cl(A) = A/\Tl(A)$,
which does not have any $\ell$-torsion.
Let $\Tl \cK_m$ be the torsion subsheaf of the Zariski sheaf $\cK_m$, which is the Zariski sheaf associated
to the presheaf $U \mapsto \Tl(K_m(U))$, and let $\Cl \cK_m$ be the Zariski sheaf associated to the presheaf
$U \mapsto \Cl(K_m(U))$. Then we have an exact sequence
$$
0 \rightarrow \Tl \cK_m \rightarrow \cK_m \rightarrow \Cl \cK_m \rightarrow 0 \,,\leqno(5.3.1)
$$
where $\Cl \cK_m$ is an $\ell$-torsion-free sheaf. By the last claim of Theorem 4.3, the functors sending $X$ to $H^i(X,\Tl \cK_m)$
or to $H^i(X,\Cl \cK_m)$ are sufficiently rigid for smooth varieties, and the sequence (5.3.1) induces a long exact cohomology sequence
$$
\ldots \rightarrow H^i(X,\Tl \cK_m) \rightarrow H^i(X,\cK_m) \rightarrow
H^i(X,\Cl \cK_m) \rightarrow H^{i+1}(X,\Tl \cK_m) \rightarrow \ldots \,,\leqno(5.3.2)
$$
which is functorial in $X$ for the functorialities of a sufficiently rigid functor, and gives rise
to an exact sequence of sufficiently rigid functors in an obvious way.

\medskip\noi
On the other hand we can apply the functors $\Tl$ and $\Cl$ to the complexes $E^{\ast,q}_1(X)$
arising from the spectral sequence (4.3.1) and get an exact sequence of complexes
$$
0 \rightarrow \Tl(E^1_{\ast,q}(X)) \rightarrow E^1_{\ast,q}(X) \rightarrow \Cl(E^1_{\ast,q}(X)) \rightarrow 0 \,. \leqno(5.3.3)
$$
whose homology gives the groups $E^2_{p,q}(X)_\Tl$, $E^2_{p,q}(X)$ and $E^2_{p,q}(X)_{\Cl}$, respectively.

\begin{lem}\label{Lemma 5.4}
Let $F$ be any functor from abelian groups to abelian groups. Then the natural maps induce isomorphisms
for all $k$-varieties $X$ and all natural numbers $n$
$$
\begin{array}{rcl}
E^2_{p,q}(X)_F[\ell^n] \mathop{\longrightarrow}\limits^{\alpha^\ast} E^2_{p,q}(X_K)_F [\ell^n] & and &
E^2_{p,q}(X)_F/\ell^n \mathop{\longrightarrow}\limits^{\alpha^\ast}  E^2_{p,q}(X_K)_F/\ell^n
\end{array}\,,
$$
where $E^2_{p,q}(Y) = H_p(FE^1_{\ast,q}(Y))$ (see (5.1.8).
\end{lem}

\bigskip\noindent\textbf{Proof}
If $Z \subset X$ is a closed subscheme and $U = X\setminus Z$ is the closed complement,
then one has exact sequences and a morphism of exact sequences
$$
\xymatrix{
\ldots\ar[r] & E^2_{p,q}(Z_K)_F \ar[r] & E^2_{p,q}(X_K)_F\ar[r] &
E^2_{p,q}(U_K)_F \ar[r] & E^2_{p-1,q}(Z_K)_F \ar[r] & \ldots\\
\ldots\ar[r] & E^2_{p,q}(Z)_F\ar[r] \ar[u]_{\alpha^\ast}& E^2_{p,q}(X)_F\ar[r]\ar[u]_{\alpha^\ast} & E^2_{p,q}(U)_F\ar[r]\ar[u]_{\alpha^\ast} &
E^2_{p-1,q}(Z)_F\ar[r]\ar[u]_{\alpha^\ast} & \ldots \,,
}\leqno(5.4.1)
$$

\noi
With this we proceed by induction on dimension.
First assume that $k$ is perfect. Then the claim is true for all $k$-varieties $X$ of dimension $0$, because these
are smooth and we have isomorphisms $E_2^{0,q}(X)=K_{-q}(X)=H^0(X,\cK_{-q})$.
If now $X$ is arbitrary, then, by perfectness of $k$, there is an smooth open dense subvariety $U \subseteq X$.
Here, by Theorem 4.3 and with the notation as in the proof of Theorem 4.3 we have a commutative diagram
\begin{footnotesize}
$$
\begin{array}{rl}
{\xymatrix{
E_2^{p,q}(U_K)_F \ar[r]^{\cong}          &  H^p(U_K,F\cK_{-q})    \\
E_2^{p,q}(U)_F \ar[r]^{\cong} \ar[u]   &  H^p(U, F\cK_{-q})\ar[u]
}}
\end{array} \leqno(4.5.5)
$$
\end{footnotesize}

\noi
with horizontal isomorphisms. By assumption the right vertical map is an isomorphism.
Therefore the left vertical maps are isomorphism as well.
Since the dimension of $Z$ is strictly smaller than that of $X$, we can assume that the claim holds for
$Z$, and hence we get the claim for $X$ by the diagram (4.5.4). If $k$ is not perfect, we can pass to the perfect hull $k^{per}$
of $k$, and a trace argument shows that the groups are the same for $k$ and $k^{per}$ and for $K$ and the composition $Kk^{per}$
after inverting the characteristic $p$ of $k$.

\bigskip\noi
Now, from the exact sequence (4.5.3) of complexes for $X$ and $X_K$ we get a commutative diagram with exact rows
$$
\xymatrix{
\ldots  E^{p,q}_2(X_K)_\Tl \ar[r]   &  E^{p,q}_2(X_K) \ar[r]  & E^{p,q}_2(X_K)_\Cl \ar[r] & E^{p+1,q}_2(X_K)_\Tl  \ldots\\
\ldots  E^{p,q}_2(X)_\Tl \ar[r] \ar[u]_{\alpha^\ast} &  E^{p,q}_2(X) \ar[r]\ar[u]_{\alpha^\ast} & E^{p,q}_2(X)_\Cl  \ar[r]\ar[u]_{\alpha^\ast} &
E^{p+1,q}_2(X)_\Tl \ar[u]_{\alpha^\ast} \ldots \,,
}\leqno(4.5.6)
$$
in which the groups $E^{\ast,\ast}_2(X)_\Tl$ and $E^{p,q}_2(X_K)_\Tl$ are torsion groups. This induces
a commutative diagram with exact rows
$$
\xymatrix{
\ldots  E^{p,q}_2(X_K)_\Tl \ar[r]   & \Tl E^{p,q}_2(X_K) \ar[r]  & \Tl (E^{p,q}_2(X_K)_{\Cl}) \ar[r] & E^{p+1,q}_2(X_K)_\Tl  \ldots\\
\ldots  E^{p,q}_2(X)_\Tl \ar[r] \ar[u]_{\alpha^\ast}^{(1)} & \Tl E^{p,q}_2(X) \ar[r]\ar[u]_{\alpha^\ast}^{(2)} & \Tl(E^{p,q}_2(X)_\Cl)  \ar[r]\ar[u]_{\alpha^\ast}^{(3)} &
E^{p+1,q}_2(X)_\Tl \ar[u]_{\alpha^\ast}^{(4)} \ldots \,,
}\leqno(4.5.7)
$$
By Claim 4.6 the maps $\alpha^\ast$ are isomorphisms at the positions (1), (3) and (4),
and the 4-Lemma implies that $\alpha^\ast$ is surjective at position (2). Therefore it
suffices to show that this map is injective. For this it suffices to show that the map
$$
\alpha^\ast: E^{p,q}_2(X)[\ell] \mathop{\longrightarrow}\limits^{\alpha^\ast} E^{p,q}_2(X_K)[\ell]
$$
is injective. For this we apply

\medskip\noi The following observation will be useful in the next section.

\begin{lem}\label{Lemma 4.6} Let $V$ be a sufficiently rigid functor
on the category $Sch^{noeth}/k$, where $k$ is a perfect field.
Then properties 4.2 (i) and (ii) hold if and only if they hold for
rings $A$ which are localizations of smooth $k$-schemes.
\end{lem}

\medskip\noindent \textbf{Proof} (cf. \cite{Qui} \S 7, Proof of thm. 5.11)\quad
Let $L$ be an extension of $k,Y= Spec ~R$ a smooth
affine $L$-scheme, and $A$ a semi-local ring obtained by
localizing $R$ in a finite set of primes $S$. Then there exists a
subfield $L'$ of $L$ finitely generated over $k$, a smooth
$L'$-Algebra $R'$ and a finite subset $S'$ of $Spec~R'$ such that
$R=L\otimes_{L'} R'$ and such that the primes in $S$ are the base
extensions of the primes in $S'$. This follows by applying
\cite{EGA IV} 8.8.2, 8.10.5, 8.7.3 and 17.7.8 to the family $(L_i, i\in I)$
of subfields of $L$ which are finitely generated over $k$. Let
$L_j~~~(j\in J)$ be the subfamily of those fields which contain
$L'$, and for each $j\in J$ let $S_j$ be the set of primes in
$R_j=L_j\otimes_{L'}R'$ obtained by tensoring the primes in $S'$
with $L_j$, and let $A_j$ be the localization of $R_j$ in $S_j$.
Then $A=\mathop{\lim}\limits_{\longrightarrow\,j} A_j$, so by the
limit property 1.6 (e) for $V$ and the $H^\nu(-,\mathcal V)$, we
see that it suffices to show 4.2 (i) and (ii) for the localization
$A$ of a smooth $L$-scheme $Y$, where $L$ is finitely generated
over $k$. But since $k$ is perfect, every such $L$ is the function
field of a smooth $k$-scheme; hence $A$ is the localization of a
smooth $k$-scheme, and the claim follows.

\section{Poincar\'e duality theories and $\cH$-cohomology}

To show the properties 4.1 (i) and (ii) for $\mathcal H$-cohomology,
where $\mathcal H$ is the Zariski sheaf associated
to \'etale cohomology, we will consider more generally the case
of a twisted Poincar\'e duality theory (with supports) as introduced
by Bloch and Ogus \cite{BO} (1.3). It encodes the usual properties of a
cohomology theory with supports for algebraic schemes over a field
$k$, with an associated (Borel-Moore type) homology theory. We
recall the axioms, since we need to consider them more closely.

\begin{defi}\label{Definition5.1} Let $k$ be a field, and
let $\cS$ be a category of algebraic $k$-schemes such that
$Y\in~Ob(\cS)$ if $X\in~Ob(\cS)$ and $Y\subset X$ is locally
closed.

\bigskip\noindent (1) Let $\cS^\ast$ be the category whose objects
are the closed immersions $Y\subset X$ in $\cS$ and whose
morphisms are cartesian squares
$$
\begin{array}{cccc} Y' &\subset& X'\\ \downarrow& &\downarrow\\
                    Y  &\subset& X & .
\end{array}
$$
A twisted cohomology theory (with supports) is a sequence (indexed
by $j\in\mathbb Z$) of contravariant functors
$$
\begin{array}{rcl}
\cS \hspace{0,6cm} & \longrightarrow & \mbox{(graded abelian groups)}\\
(Y\subset X)       & \mapsto         & \hspace{0,8cm}\oplus_i H_Y^i (X,j)
\end{array}
$$
satisfying the following properties

\v\noi (a) (long exact cohomology sequence) For $Z\subset Y\subset X$, there is a long exact sequence
$$
\ldots\rightarrow H^i_Z(X,j)\rightarrow H^i_Y(X,j)\rightarrow
H^i_{Y-Z}(X\smallsetminus Z,j)\rightarrow
H^{i+1}_Z(X,j)\rightarrow\ldots,
$$
functorial with respect to morphisms
$$
\begin{array}{ccccc}
Z'& \subset & Y' & \subset & X' \\
\downarrow & & \downarrow & & \downarrow\\
Z & \subset & Y & \subset & X
\end{array}
$$
in the obvious way.

\v\noi (b) (excision) If $Z\subset X\in~~Ob~\mathcal V^\ast$ and if
$U\subset X$ is open in $X$ and contains $Z$, then the map
$H^i_Z(X,j)\rightarrow H^i_Z(U,j)$ is an isomorphism.

\bigskip\noindent (2) Let $\cS_\ast$ be the category with
$~Ob(\cS_\ast) = ~Ob\cS$ but whose arrows are only the proper
morphisms in $\cS$. A twisted homology theory is a sequence
(indexed by $b\in\mathbb Z$ of covariant functors
$$
\begin{array}{rcl}
\varphi_\ast & \rightarrow & \hbox{(graded abelian groups)}\\
X & \shortmid\hspace{-0,09cm}\rightsquigarrow & \hspace{0,7cm} \scriptstyle\bigoplus\limits_a \textstyle H_a(X,b)\\
\end{array}
$$
such that the following properties hold.

\v\noi (c) If $\alpha:X'\rightarrow X$ is \'etale, there is a map
$$
\alpha^\ast: H_a(X,b)\rightarrow H_a(X',b)~~,
$$
such that $(\alpha~\alpha')^\ast=\alpha'^\ast \alpha^\ast$ for
$\alpha': X''\rightarrow X'$ \'etale.

\v\noi (d) If the diagram below on the left is cartesian, with
proper $f$ and $g$, and \'etale $\alpha$ and $\beta$, then the
diagram on the right commutes.
$$
\xymatrix{ X'\ar[d]_g\ar[r]^\beta& X\ar[d]^f & &
H_a(X,b)\ar[d]_{f_\ast} \ar[r]^{\beta^\ast} &
H_a(X',b)\ar[d]^{g_\ast}\\
Y' \ar[r]^\alpha & Y & & H_a(Y,b)\ar[r]^{\alpha^\ast} & H_a(Y',b)&
.}
$$

\v\noi (e) If $i:Z\hookrightarrow X$ is a closed immersion in
$\mathcal V$, with open complement $\alpha: V\hookrightarrow X$,
then there are long exact sequences
$$
\ldots\rightarrow H_a(Z,b)\stackrel{i_\ast}{\rightarrow}
H_a(X,b)\stackrel{\alpha^\ast}{\rightarrow} H_a(U,b)\rightarrow
H_{a-1}(Z,b)\rightarrow\ldots~~~.
$$
which satisfy the following compatibilities (NB only the
commutativity of the squares (1) and (2) below is a new
statement):

\v\noi (f) If $f:X'\rightarrow X$ is a proper morphism,
restricting to $f':Z'\rightarrow Z$ for a closed subscheme
$Z'\subseteq X$ then the diagram
$$
\xymatrix{ \ldots\ar[r] & H_a(Z',b) \ar[r]\ar[d]^{f'_\ast} &
H_a(X',b)\ar[r]\ar[d]^{f_\ast} & H_a(U',b)
\ar[r]\ar[d]^{f_\ast\alpha^\ast~~\displaystyle (1)} &
H_{a-1}(Z',b)
\ar[r]\ar[d]^{f'_\ast} & \ldots\\
\ldots\ar[r] & H_a(Z,b)\ar[r] & H_a(X,b)\ar[r] & H_a(U,b)\ar[r] &
H_{a-1}(Z,b)\ar[r] & \ldots}
$$
commutes, where $\alpha:f^{-1}(U)=X'-f^{-1}(Z)\rightarrow
X'-Z'=U'$ is the open immersion.

\v\noi (g) If $\alpha:X'\rightarrow X$ is \'etale, then the
diagram
$$
\xymatrix{ \ldots\ar[r] & H_a(Z',b) \ar[r] & H_a(X',b)\ar[r] &
H_a(U',b) \ar[r] & H_{a-1}(Z',b)
\ar[r] & \ldots\\
\ldots\ar[r] & H_a(Z,b)\ar[r] \ar[u]_{\alpha^\ast}&
H_a(X,b)\ar[r]\ar[u]_{\alpha^\ast} &
H_a(U,b)\ar[r]\ar[u]_{\alpha^\ast\quad\quad\displaystyle(2)} &
H_{a-1}(Z,b)\ar[r]\ar[u]_{\alpha^\ast} & \ldots & ,}
$$
commutes, where $Z'=f^{-1}(Z)$ and $U'=f^{-1}(U)~.$

\bigskip\noi (3) A Poincar\'e duality theory is given by a cohomology
and homology theory as above, together with the following
structures.

\v\noi (h) (cap product) For $Y\subset X\in~Ob~\varphi^\ast$ there
is a pairing
$$
\cap:\; H_a(X,b)\otimes H^i_Y(X,j)\rightarrow H_{a-i}(Y,b-j)~,
$$
compatible with \'etale pull-backs, in the obvious way.

\v\noi (i) (projection formula) For a cartesian diagram below on
the left, with proper $f$, the diagram on the right commutes.
$$
\xymatrix{ Y'\ar[d]\ar@{^{(}->}[r] & X'\ar[d]^f & & \cap:\; H_a(X',b)\ar[d]^{f_\ast}& \otimes &  H^i_{Y'}(X',j)\ar[r] & H_{a-i}(Y',b-j)\ar[d]^{f_\ast}\\
           Y \ar@{^{(}->}[r]       & X          & & \cap:\; H_a(X,b)                & \otimes &  H^i_{Y}(X,j)\ar[u]_{f^\ast}\ar[r]   & H_{a-i}(Y,b-j) & .}
$$

\v\noi (j) (fundamental class) If $X\in ~Ob~\varphi$ is irreducible
and of dimension $d$, then there is a canonical element $\eta_X\in
H_{2d} (X,d)$. If $\alpha:X'\rightarrow X$ is \'etale, then
$\alpha^\ast\eta_X=\eta_X'~.$

\v\noi (k) (Poincar\'e duality) If $X\in~Ob~\varphi$ is smooth of
pure dimension $d$, and $Y\subset X$ is a closed subscheme, then
cap-product with $\eta_X$ induces an isomorphism
$$
\eta_X\cap~:~H^i_Y(X,j)\stackrel{\sim}{\longrightarrow}
H_{2d-i}(Y,d-j).
$$

\v\noi (l) (Principal triviality) Let $i: W\hookrightarrow X$ be a
smooth principal divisior in the smooth scheme $X$. Then $i_\ast
\eta_W=0~.$
\end{defi}

\bigskip\noi Write $H^i(X,j)$ for $H^i_X(X,j).$ By results of Bloch and Ogus we then have:

\begin{prop}\label{Proposition 5.2} Let $k$ be a perfect field, and let
$$
(Z\subset X)\mapsto H^\ast_Z(X,\ast)~~ , ~~ X\mapsto H_{\ast}(X,\ast)
$$
be a twisted Poincar\'e duality theory on the category $Sch^{alg}/k$ of
all algebraic $k$-schemes.
%Let $\cH^i(j)$ be the Zariski sheaf associated to the presheaf $U \mapsto H^i(U,j)$.
Then for all $i,j\in\mathbb Z$ the functor $V:(Sch^{alg}/k)^0\rightarrow
Ab$ defined by
$$
V(Y)=H^i(Y,j)
$$
satisfies property 4.1 (i) for semi-local regular rings $A$ which
are localizations of smooth $k$-schemes.
\end{prop}

\noi Since $V$ is a priori only defined on $k$-schemes of finite
type, this statement has to be interpreted in the following way:
the functors $V$ and $H^\nu(-,\mathcal V)$ are defined on $Spec~A$
by taking the limit of the value groups at all opens $U\subset X$
containing $Spec~A$ (or, equivalently, the maximal ideals of $A$).

\medskip\noindent \textbf{Proof of 5.2} For any algebraic $k$-scheme
$X$ denote by $X_{(p)}$ the set of points $x\in X$ of dimension
$p$, and for $x\in X$ put
$$
H_a(x,b)=\mathop{\longrightarrow}\limits^{\lim}_{U\subset Z~{\rm open}}
H_a(U,b)~~~,
$$
where $Z=\overline{\{x\}}$ is the Zariski closure of $x$. Then
Bloch and Ogus construct a homological spectral sequence
$$
E^1_{p,q}=E^1_{p,q}(X,b)=\mathop{\oplus}\limits_{x\in X_{(p)}}
H_{p+q}(x,b)~\Rightarrow~ H_{p,q}(X,b)\leqno(5.2.1)
$$
as follows (\cite{BO} (3.7)). Let $Z_p=Z_p(X)$ be the set of all closed
subsets $Z\subseteq X$ of dimension $\leq p$, ordered by
inclusion, and put
$$
H_a(Z_p(X),b)=\mathop{\longrightarrow}\limits^{\lim}_{Z\in Z_p} H_a(Z,b)~~~.
$$
Then by 5.1 (e) one gets exact sequences
$$
\ldots\rightarrow H_a(Z_{p-1},b)\stackrel{i}{\rightarrow}
H_a(Z_p,b)\stackrel{j}{\rightarrow} \mathop{\oplus}\limits_{x\in
X_{(p)}} H_a(x,b)\stackrel{k}{\rightarrow}
H_{a-1}(Z_{p-1},b)\rightarrow\ldots \; .\leqno(5.2.2)
$$
The method of exact couples now gives the desired spectral
sequence. Note that by definition the differential $d^1_{p,q}$ is
the composition
$$
d^1_{p,q}:E^1_{p,q}\stackrel{k}{\rightarrow}
H_{p+q-1}(Z_{p-1},b)\stackrel{j}{\rightarrow} E^1_{p-1,q}\quad. \leqno(5.2.3)
$$
Moreover, if $X$ is of pure dimension $d$, one has a complex
$$
0\rightarrow H_a(X,b)\stackrel{\varepsilon}{\rightarrow}
E^1_{d,a-d}(X,b) \stackrel{d^1}{\rightarrow}
E^1_{d-1,a-d}(X,b)\stackrel{d^1}{\rightarrow}\ldots\quad,
\leqno(5.2.4)
$$
in which $\varepsilon$ is the edge morphism. By 5.1 (g) the
sequences (5.2.4) for all opens $U \subseteq X$ form a complex of
Zariski presheaves. Let $\mathcal H_a(b)$ and $\mathcal
E^1_{p,q}(b)$ be the Zariski sheaves associated to $U\mapsto
H_a(U,b)$ and $U\mapsto E^1_{p,q} (U,b)=
\mathop{\oplus}\limits_{x\in U_{(p)}} H_{p+q}(x,b)$, respectively,
so that we get a complex of Zariski sheaves
$$
0\rightarrow\mathcal H_a(b)\rightarrow\mathcal E^1_{d,a-d}(b)
\rightarrow\mathcal E^1_{d-1,a-d}(b)\rightarrow\ldots
\quad.\leqno(5.2.5)
$$
It is clear from the definition that $U\mapsto E^1_{p,q} (U) $ is
already a sheaf, and is flabby.

\medskip Now let $X$ be smooth. Then Bloch and Ogus \cite{BO} show that for any
finite set $S \subset X$ which is contained in an affine open, the
sequence (5.2.3) becomes exact after passing to the limit over all
opens $V\subset X$ containing $S$. In fact, this is equivalent to
the statement that all maps $i:H_a(Z_{p-1}(X),b)\rightarrow
H_a(Z_p (X),b)$ vanish after passing to the limit over such $U$,
and this is shown in \cite{BO}, section 4 and 5. (In the claim on p.191
loc. cit. only the case of a one-element set $S$ is stated, but
the proof works more generally, since the trick of Quillen quoted
loc. cit. is valid for a finite $S$ as above.) In particular,
(5.2.5) is an exact sequence of Zariski sheaves, hence a
resolution of $\mathcal H_a(b)$ by the complex with flabby
components
$$
\mathcal E^1_{d-\ast,a-d}(b)\quad.\leqno(5.2.6)
$$
As a consequence, by applying the functor $H^0(U,-)$ we get a canonical isomorphism
$$
\alpha:  H^\nu(U,\mathcal H_a(b))= H^\nu(H^0(U,\mathcal E^1_{d-\ast,a-d}(b))) = E^2_{d-\nu,a-d}(U,b)\leqno(5.2.7)
$$
for every open $U$ in $X$ and $\nu\geq 0$. Hence, if $A$ is a semi-local ring obtained by localizing $X$,
then by exactness of (5.2.4) we have $H^\nu (Spec~A,\mathcal H_a(b))=0$ for all $\nu>0$, and the map
$$
H_a(Spec~A,b)\rightarrow E^2_{a,d-a}(Spec~A,b)=H^0(Spec~A,\mathcal H_a(b))
$$
is an isomorphism. Since the presheaves $U\mapsto H_a(U,b)$ and
$U\mapsto H^{2d-a}(U,d-b)$ are isomorphic by 4.1 (h), (j) and (k),
the proposition follows.

\bigskip For the treatment of the homotopy property 4.1 (ii) for
$\mathcal H$-cohomology we need the following extended version of
a Poincar\'e duality theory.

\begin{defi}\label{Definition 5.3} Let $k$ be a field. A
twisted Poincar\'e duality theory
$$ (Z\subset X)\mapsto H^\ast_Z(X,\ast) ~~,~~ X\mapsto H_\ast(X,\ast)
$$
on $Sch^{alg}/k$ is called an extended Poincar\'{e} duality theory, if for every flat
morphism $f:X'\rightarrow X$ which is equidimensional of dimension
$m$ (i.e., whose fibres are either empty or equidimensional of
dimension $m$ , cf. \cite{EGA IV}, 13.3) there are functorially
associated maps
$$
f^\ast: H_a(X,b)\rightarrow H_{a+2m}(X',b+m)~,
$$
agreeing with the pull-back maps in 5.1 (c) for \'etale $f$ and
$m=0$, such that the following further properties hold.

\v\noi
(m) If $X$ and $X'$ are irreducible, then $f^\ast\eta_X =\eta_{X'}$.

\v\noi
(n) If $Z\subset X$ is a closed subscheme, and
$Z'=Z\times_X X' \subset X'$, then the following diagram commutes.
$$
\xymatrix{
\cap: & H_{a+2m}(X',b+m) & \otimes & H^i_{Z'}(X',j) \ar[r]              & H_{a-i+2m}(Z',b-j+m)& \\
\cap: & H_a(X,b)\ar[u]_{f^\ast} & \otimes & H^i_{Z}(X,j)\ar[r]\ar[u]_{f^\ast}  & H_{a-i}(Z,b-j)\ar[u]_{f^\ast} & .}
$$

\v\noi
(o) If $Z$ is closed in $X$ and $U$ is the open complement,
then one has a commutative diagram
\begin{small}
$$
\xymatrix{ ..\ar[r] & H_{a{+}2m}(Z',b{+}m)\ar[r] &
H_{a{+}2m}(X'{,}b{+}m)
\ar[r] & H_{a{+}2m}(U'{,}b{+}m) \ar[r] & H_{a-1+2m}(Z',b{+}m)\\
..\ar[r] & H_a(Z,b) \ar[u]_{f^\ast} \ar[r] &
H_a(X,b)\ar[u]_{f^\ast} \ar[r] & H_a(U,b)\ar[u]_{f^\ast} \ar[r] &
H_{a-1}(Z,b)\ar[u]_{f^\ast}\longrightarrow .. }
$$
\end{small}
where $Z'=Z\times_X X'$ and  $U'=f^{-1}(U)$.
\end{defi}

\medskip\noi We shall give some examples below. First we note

\begin{prop}\label{Proposition 5.4} Let $k$ be a perfect
field, and let $(Z\subset X)\mapsto H^\ast_Z (X,\ast),X\mapsto
H_.(X,\ast)$ be an extended Poincar\'e duality theory on
$Sch^{alg}/k$. Assume that the following ``homotopy invariance''
holds:

\v\noi (p) For every smooth $k$-scheme $X$ the maps
$$
p^\ast:H^i(X,j)\longrightarrow H^i(\mathbb A^1_X,j)
$$
induced by the projection $p:\mathbb A^1_X\rightarrow X$ are
isomorphisms for all $i,j\in\mathbb Z$.

\v\noi Then for all $i,j\in\mathbb Z$ the contravariant functor
$V$ on $Sch^{alg}/k$ given by
$$
V(Y)~~=~~H^i(Y,j)
$$
satisfies property 4.1 (ii) for local rings $A$ which are local
rings $\mathcal O_{Y,y}$ of smooth $k$-schemes $Y$.
\end{prop}

\medskip\noindent \textbf{Proof} (The values on $Spec~A$ are
defined as in Proposition 5.2). Let $f:X'\rightarrow X$ be a flat
morphism in $Sch^{alg}/k$ which is equidimensional of dimension
$m$.

\medskip We first note that by 5.3 (o) the flat pull-backs induce maps
between the sequences (5.2.2) for $X$ and $X'$, with appropriate
shift of degrees (for $Z\subset X$ of dimension $p$ the preimage
$f^{-1}(Z)=Z\times_X X'$ is of dimension $p+m$). This induces a
map of spectral sequences from (5.2.1) (with the indicated shift of degrees)
$$
\xymatrix{
E^1_{p+m,q+m}(X',b+m) \ar@{=>}[r] & H_{p+q+zm}(X',b+m)\\
E^1_{p,q}(X,b)\ar[u]_{f^\ast}\ar@{=>}[r] &
H_{p+q}(X,b)\ar[u]_{f^\ast} & .}\leqno(5.4.1)
$$
In particular, there are natural pull-back maps
$$
f^\ast:E^2_{p,q}(X,b)\rightarrow E^2_{p+m,q+m}(X',b+m)~~~~.\leqno(5.4.2)
$$

\medskip On the other hand, if $X$ and $X'$ are irreducible of
dimensions $d$ and $d'$, respectively, (so that $m=d'-d)$, then by
5.3 (m) and (n) the diagram
$$
\xymatrix{
H^i(X',j)\ar[r]^{\hspace{-0,7cm}\eta_{X'}\cap} & H_{2d'-i}(X',d'-j)\\
H^i(X,j)\ar[u]^{f^\ast}\ar[r]^{\hspace{-0,5cm}\eta_X \cap} &
H_{2d-i}(X,d-j)\ar[u]_{f^\ast}}\leqno(5.4.3)
$$
is commutative. Now let $X$ and $X'$ be smooth. Then the
horizontal maps in the above diagram are isomorphisms by
Poincar\'e duality 3.1 (k). The same is true for open subschemes;
hence the pull-back map $f^\ast: H^\nu(X,\mathcal
H^i(j))\rightarrow H^\nu(X',\mathcal H^i(j))$ can be identified
with a pull-back map
$$
f^\ast:H^\nu(X,\mathcal H_{2d-i}(d-j))\rightarrow H^\nu(X',\mathcal H_{2d'-i}(d'-j))
\leqno(5.4.4)
$$
which is defined in a way analogous to (4.2.2), using the
compatibility of flat pull-backs with open immersions. Here
$\mathcal H_a(b)$ is a Zariski sheaf on $X$ or $X'$ associated to
the presheaf $U\mapsto H_a(U, b).$

\medskip Furthermore, it follows from (5.4.1) that, via the
isomorphisms (5.2.6) for $X$ and $X'$, the pull-back map (5.4.4)
can be identified with the pull-back map (5.4.2) for
$(p,q,b,m)=(d-\nu,d-i,d-j, d'-d)$. Thus, in view of Remark 4.2 (b),
the Proposition follows from part (ii) of the Lemma below,
applied to smooth $k$-schemes $X$.

\begin{lem}\label{Lemma 5.5} In the situation of Proposition
5.4, let $X$ be any algebraic $k$-scheme (not necessarily smooth),
and let $p:\mathbb A^1_X\rightarrow X$ be the affine line over
$X$. Then the following holds.

\v\noi (i) The flat pull-back maps $p^\ast:H_a(X,b)\rightarrow
H_{a+2}(\mathbb A^1_X,b+1)$ are isomorphisms for all
$a,b\in\mathbb Z$.

\v\noi (ii) The pull-back maps $p^\ast: E^2_{p,q}(X,\,
b)\rightarrow E^2_{p+1,\, q+1}(\mathbb A^1_X,\, b+1)$ are
isomorphisms for all $p,q,b\in\mathbb Z$.
\end{lem}

\bigskip\noindent\textbf{Proof} We proceed by induction on $dim(X)$, and we
may and will consider only reduced schemes (since $H_a(X,\, b) =
H_a(X_{red},\, b)$ by 4.1 (e)). For dim$(X)=0$ we then may
assume that $X=Spec~ K$ for a finite extension $K$ of $k$,
necessarily separable since $k$ is perfect. Then (5.4.1) gives a
commutative diagram with exact top row
$$
\xymatrix{ 0 \ar[r] & E^2_{0,a}(\mathbb A^1_K,b+1) \ar[r] &
H_{a+2}(\mathbb A^1_K,b+1) \ar[r] & E^2_{1,a-1}(\mathbb A^1_K,b+1) \ar[r] & 0\\
&& H_a(Spec(K),b) \ar[u]_{p^\ast} \ar[r]^\sim & E^2_{
0,a}(Spec(K), b) \ar[u]_{p^\ast}} \leqno(5.5.1)
$$
in which the middle vertical map is an isomorphism by the
assumption 5.4 (p) and Poincar\'e duality (cf. (5.4.3)). On the
other hand, the right hand vertical map is injective: it can be
identified with
$$
p^\ast: H^0(Spec~(K),\mathcal H^{-a}(-b))\rightarrow H^0(\mathbb
A^1_K,\mathcal H^{-a}(-b))~~,
$$
and any $K$-rational point $s: Spec~K\rightarrow\mathbb A^1_K$
gives a left inverse $s^\ast$ of $p^\ast$. Putting this together,
we deduce that both vertical maps are isomorphisms and that
$E^2_{0,a}(\mathbb A^1_K,b+1)=0~.$ This shows (i) and (ii) in this
case, since $E^2_{p,q}(Spec(K),b)=0$ for $p\not= 0$ and
$E^2_{p,q}(\mathbb A^1_K,b)=0$ for $p\not= 0\,, -1$.

\medskip If $X$ has positive dimension (and is reduced), then there
is a dense open subscheme $U\subseteq X$ which is smooth (because
$k$ is perfect). Since (i) holds for $U$ by Poincar\'e duality and
assumption, and for $Z=X-U$ by induction hypothesis, it holds for
$X$ by 5.3 (o) and the five-lemma.

\medskip For (ii) we observe the following: if
$Z\mathop{\hookrightarrow}\limits^i X$ is a closed subscheme and
$U=X-Z\stackrel{j}{\hookrightarrow} X $ is the open complement,
then one has an exact sequence of complexes
$$
0\rightarrow E^1_{\cdot,q}
(Z,b)\stackrel{i_\ast}{\rightarrow}E^1_{\cdot,q}
(X,b)\stackrel{j^\ast}{\rightarrow}E^1_{\cdot,q} (U,b)\rightarrow
0.\leqno(5.5.2)
$$
Here $j^\ast$ comes from (5.4.1) for $f=j$, and $i_\ast$ comes
from the contravariance of the spectral sequence (5.2.1) for the
proper morphism $i$, which is an immediate consequence of 3.1 (f).
The exactness of the sequence follows from the easily verified
fact that in degree $p$ the sequence is given by
$$
0\rightarrow
\mathop{\oplus}\limits_{x\in{Z_{(p)}}}H_{p+q}(\kappa(x),b)
\stackrel{i_\ast}{\longrightarrow}\mathop{\oplus}\limits_{x\in{X_{(p)}}}
H_{p+q}(\kappa(x),b)\stackrel{j^\ast}{\longrightarrow}
\mathop{\oplus}\limits_{x\in{U_{(p)}}}H_{p+q}(\kappa(x),b)\longrightarrow
0~,
$$
where $i_\ast$ and $j^\ast$ are the obvious inclusion and
projection, respectively. There is a corresponding exact sequence
for the triple $\mathbb A^1_Z\hookrightarrow \mathbb
A^1_X\hookleftarrow\mathbb A^1_U$, and both exact sequences are
connected by the pull-back maps for the projections $p_X:\mathbb
A^1_X\rightarrow X~,~ p_U$ and $p_Z$, in a commutative map.
Passing to the cohomology, we obtain a commutative diagram with
long exact rows
\begin{small}
$$
\xymatrix{..\ar[r] E^2_{p+1,q+1}(\mathbb A^1_Z,b\!+\!1) \ar[r] & E^2_{p+1,q+1}(\mathbb A^1_X,b\!+\!1)\ar[r]
& E^2_{p+1,q+1}(\mathbb A^1_U,b\!+\!1)\ar[r] & E^2_{p,q+1}(\mathbb A^1_Z,b\!+\!1) ..\\
..\ar[r] E^2_{p,q}(Z,b)\ar[r] \ar[u]_{p_Z^\ast}& E^2_{p,q}(X,b)\ar[r] \ar[u]_{p_X^\ast}&E^2_{p,q}(U,b)\ar[r]\ar[u]_{p_U^\ast} &
E^2_{p-1,q}(Z,b)\rightarrow\ar[u]_{p_Z^\ast}..}\leqno(5.5.3)
$$
\end{small}

\v\noi By induction hypothesis we may assume that all $p^\ast_Z$ are
isomorphisms for $\dim Z<\dim$ $X$. Hence it suffices to show that
the $p^\ast_U$ become isomorphisms after passing to the limit over
all dense opens $U\subset X$. Moreover, we may restrict to the
case that we consider the limit over all open neighbourhoods of a
fixed generic point $\eta$. If $K=\kappa(\eta)$ is the function
field of the corresponding connected component of $X$, then we
obtain formally the same diagram as (5.5.1), by putting
$H_a(K,b)= \mathop{\lim}\limits_{\longrightarrow} H_a(U,b), \,
H_a(\mathbb A^1_K,b) = \mathop{\lim}\limits_{\longrightarrow} H_a(\mathbb A^1_U,b)$
etc. By reasoning in a completely similar way as in the case of
a finite extension $K$ of $k$, we deduce that the pull-back map
$$
p^\ast:E^2_{p,q}(K,b)\rightarrow E^2_{p+1,q+1}(\mathbb A^1_K,b+1)
$$
is an isomorphism as wanted (the map is injective since
$p^\ast:E^2_{p,q} (U,b)\rightarrow E^2_{p+1,q+1}(\mathbb
A^1_U,b+1)$ is injective for each smooth open $U\subseteq X$).

\begin{rem}\label{Remark 5.6} The proof of 5.5 follows very
much Gillet's proof of the corresponding statement for the Chow
groups $CH_{r,s}~(X)$ (\cite{Gi} Thm. 8.3), except that the proof for
$X= Spec~(K)$ via (5.5.1) is more direct than the recursion to the
projective bundle theorem in \cite{Gi} Lemma 8.4.
\end{rem}

\noi We want to apply the above to the following example.
%For other examples of extended Poincar\'e duality theories with
%homotopy invariance, we refer to the appendix.

\begin{prop}\label{Proposition 5.7} Let $k$ be a field, and
fix $n\in\mathbb N$ invertible in $k$ The following functors form
an extended Poincar\'e duality theory on $Sch^{alg}/k$, and the
properties 4.1 (i) and (ii) hold for them (In particular, by
Poincar\'e duality, the homotopy invariance
 5.4 (p) holds for them):
$$
\begin{array}{lcll}
 H^i_Z(X,j) & = & H^i_Z(X,\mathbb Z/n(j)) & \hbox{(\'etale
cohomology with supports)},\\ H_a(X,b)& = & H_a(X,\mathbb
Z/n(b)):= H^{-a}(X,a^!_X \mathbb Z/n(-b))& \hbox{(\'etale
homology)}
\end{array}
$$
Here $\mathbb Z/n(j)=\mu_n^{\otimes j}$ as in 2.1,
$a_X:X\rightarrow Spec~(k)$ is the structural morphism, and
$a_X^!:D^b_c (Spec (k),\mathbb Z/n)\rightarrow D^b_c (X,\mathbb
Z/n)$ (= derived category of bounded complexes of \'etale $\mathbb
Z/n$-sheaves on $ X$ with constructible cohomology) is the right
adjoint of $R(a_X)_!$ constructed in Grothendieck-Verdier duality
\cite{SGA 4.3} XVII, 3.
\end{prop}

\medskip\noindent \textbf{Proof} That \'etale cohomology and homology
form a twisted Poincar\'e duality theory, follows from the results
in \cite{SGA 4.3}, cf. the sketch in [BO] 2.1. We now describe flat
pullbacks in the homology, for a flat morphism $f: X'\rightarrow
X$ which is equidimensional of dimension $m$ (cf. the discussion
for an algebraically closed field $k$ in Laumon's article \cite{Lau} VIII, 5).
By \cite{SGA 4.3} XVIII 2.9, there is a canonical trace morphism
$$
Tr_f:R^{2m}f_!f^\ast \mathcal F(m)\rightarrow\mathcal F~~,\leqno(5.7.1)
$$
for any \'etale $\mathbb Z/n$-sheaf $\mathcal F$ on $X$,
coinciding with the trace morphism used in 2.1 for finite $f$ (in
which case $m=0$ and $f_!=f_\ast$). Since $R^i f_!\mathcal F=0$
for $i> 2m$, this trace morphism can be regarded as a morphism
$$
Tr_f:~~ R f_!f^\ast \mathcal F (m)[2m]\rightarrow\mathcal F
$$
in $D^b_c(X,\mathbb Z/n).$ This can be extended to arbitrary
complexes $\mathcal L$ in $D^b_c(X,\mathbb Z/n)$ as follows. If
let $\mathcal F=\mathbb Z/n$ and tensor with $\mathcal L$, then by
the ``projection formula'' isomorphism
$$
\gamma_f:~~R f_!\mathcal K \otimes^L_{\mathbb Z/n} f^\ast\mathcal
L) \stackrel{\sim}{\rightarrow} R f_! (\mathcal
K\otimes^L_{\mathbb Z/n} f^\ast\mathcal L) \leqno(5.7.2)
$$
(\cite{SGA 4.3} XVII 5.2.9), we obtain a trace morphism
$$
Tr_f:R f_!f^\ast \mathcal L(m)[zm] \rightarrow\mathcal
L\leqno(5.7.3)
$$
(\cite{SGA 4.3} XVIII 2.13.2). By adjunction between $R f_!$ and
$f!$, we now get a morphism
$$
t_f:f^\ast\mathcal L(m)[2m]\rightarrow f^!\mathcal L.
$$
(\cite{SGA 4.3} 3.2.3). Applied to $\mathcal L=a^!_X \mathbb
Z/n(-b)$, for which $f^! \mathcal L=a^!_{X'}\mathbb Z/n(-b)$, this
induces the wanted pull-back maps
\begin{small}
$$
\xymatrix{ f^\ast: H^{-a}(X,a_X^!\mathbb Z/n(-b))
\ar[r]\ar@{=}[d] & H^{-a}(X',f^\ast a_X^!\mathbb Z/n(-b))
\ar[r]^{\hspace{-0,5cm}t_f} &
H^{-a-2m}(X',a^!_X,\mathbb Z/n(-b-m))\ar@{=}[d]\\
 H_a(X,\mathbb Z/n(b))& & H_{a+2m}(X',\mathbb Z/n(b+m)),}
$$
\end{small}

\noi where the first arrow is the restriction morphism which exists for
any complex of sheaves $\mathcal K$ (by the composition
$H^\ast(X,\mathcal K) \stackrel{ad_f}{\rightarrow} H^ast(X,
R f_\ast f^\ast\mathcal K)=H^\ast(X',f^\ast \mathcal K)$)

\medskip If $f$ is \'etale, then $f^!$ is identified with
$f^\ast$ via $t_f$, and the pull-back by definition is the one
used for \'etale morphisms in homology (cf. [BO] 2.1). The
functoriality of flat pull-backs is a direct consequence of the
``transitivity'' of the trace maps (3.7.1) (\cite{SGA 4.3} XVIII 2.9,
(Var 3)). This in turn implies 5.3 (m), because $\eta_X$ is the
image of $1\in \mathbb Z/n=H_0(Spec(k),\mathbb Z/n)$ under
$$
a^\ast_X: H_0(Spec(k),\mathbb Z/n)\rightarrow H_{2d}(X,\mathbb Z/n(d)~~,
$$
if $X$ is irreducible of dimension $d$ (cf. \cite{SGA 4 1/2}),
[cycle], 2.3).

\medskip For 5.3 $(n)$ we recall
that the cap product is induced by a pairing $a^!_X\mathbb
Z/n(r)\otimes a^\ast_X\mathbb Z/n(s)\rightarrow a^!_X\mathbb
Z/n(r+s)$ (cf. [BO] 2.1) which is a special case of the following,
more general one. For any morphism $g:X\rightarrow Y$ in
$Sch^{alg}/k$ and any $\mathbb Z/n$-sheaves $\mathcal F,\mathcal
G$ on $Y$ (in fact, any objects $\mathcal F,\mathcal G$ in
$D^b_c(Y,\mathbb Z/n))$, one has a pairing
$$
\varphi_g:g^!\mathcal F \otimes^L g^\ast\mathcal G\rightarrow g^!
(\mathcal F\otimes^L\mathcal G)~~,
$$
which by adjunction corresponds to the horizontal morphism making
the diagram
$$
\xymatrix{
R~g_!(g^!\mathcal F\otimes^L g^\ast \mathcal
G)\ar[r]\ar[d]_{\gamma_g} & \mathcal
F\otimes^L\mathcal G\\
R~g_! g^!\mathcal F\otimes^L\mathcal G\ar[ur]_{Ad_g\otimes id}}
$$
commutative. Here the vertical isomorphism is the projection
formula isomorphism (5.7.2), and $Ad_g:Rg_! g^!\mathcal
F\rightarrow\mathcal F$ is the adjunction map.

\medskip Now let
$f:X'\rightarrow X$ be a flat morphism which is equidimensional of
dimension $m$, and put $g'=gf:X'\rightarrow Y.$ Then we claim that
the diagram
$$
\xymatrix{ f^!g^!\mathcal F\otimes^L f^\ast g^\ast \mathcal G & =
& g'^!\mathcal F\otimes^L g'^\ast \mathcal G
\ar[r]^{\varphi_{g'}} & g'^!(\mathcal F\otimes^L \mathcal G)\\
f^\ast g^!\mathcal F\{m\}\otimes^L f^\ast g^\ast\mathcal G
\ar[u]^{t_f\otimes id}& = & f^\ast(g^!\mathcal F\otimes^L g^\ast
\mathcal G)\{m\}\ar[r]^{f^\ast \varphi_g} & f^\ast g^!(\mathcal
F\otimes^L\mathcal G)\{m\}\ar[u]_{t_f}}\leqno{(5.7.4)}
$$
commutes, where we put $\mathcal H\{m\}=\mathcal H(m)[2m]$ for a
complex $\mathcal H$. By adjunction, this amounts to the
commutating of the following two diagrams (where we have written
$f_!$ for $R f_!$, etc.)
$$
\xymatrix{ & g'_! g'^!(\mathcal F\otimes\mathcal
G)\ar@{}[d]^{\hspace{-0,4cm}\textstyle(1)}\ar[dr]^{Ad_{g'}} & & \\
g'_!(g'^!\mathcal F\otimes^L g'^\ast \mathcal
G)\ar[ur]^{g'_!\varphi_{g'}} & g'_! g'^! \mathcal
F\otimes^L\mathcal G \ar[r]^{Ad_{g'}\otimes id}
\ar[l]_{\gamma_{g'}}\ar[dr]^{g_!Ad_f\otimes id} & \mathcal F\otimes\mathcal G\\
g'_!(f^\ast g^!\mathcal F\{m\}\otimes g'^\ast \mathcal
G)\ar[u]^{g'_!(t_f\otimes id)} & g'_!f^\ast g^!\mathcal
F\{m\}\otimes\mathcal G \ar[r]^{\hspace{3mm}g_! Tr_f\otimes
id}\ar[l]_{\gamma_{g'}}\ar[u]^{g'_! t_f\otimes
id}_{\hspace{0,5cm}\textstyle(3)}\ar[d]^{\gamma_g} & g_!g^!
\mathcal F\otimes \mathcal G\ar[d]^{\gamma_g}\ar[u]_{Ad_g\otimes id} \\
g!f_!(f^\ast g^!\mathcal F\{m\}\otimes f^\ast g^\ast \mathcal
G)\ar@{=}[u]_{\hspace{2cm}\textstyle(4)} & g_!(f_!f^\ast
g^!\mathcal F\{m\}\otimes g^\ast \mathcal
G)\ar[l]_{\hspace{3mm}g_!\gamma_f}\ar[r]^{\hspace{8mm}g_!
(Tr_f\otimes id)} &
g_!(g^!\mathcal F\otimes g^\ast\mathcal G)\\
& g_!f_!f^\ast(g^!\mathcal F\otimes g^\ast\mathcal
G)\{m\}\ar@{=}[ul]\ar[ur]_{g_!Tr_f}\ar@{}[u]_{\hspace{-0,4cm}\textstyle(5)}}
$$

\vspace{1cm}

$$
\xymatrix{ g'_!g'^!(\mathcal F\otimes\mathcal
G)\ar[r]^{Ad_{g'}}\ar[dr]^{g_!Ad_f} & \mathcal F\otimes\mathcal G \\
g'_! f^\ast g^!(\mathcal F\otimes\mathcal G)\{m\}\ar[r]^{g_!
Tr_f}\ar[u]^{g'_!t_f}_{\hspace{0,5cm}\textstyle (3)} &
g_!g^!(\mathcal F\otimes\mathcal
G)\hspace{1cm}(1)\ar[u]_{Ad_g{\hspace{-1,5cm}\textstyle(2)}}
& ~g_!g^!\mathcal F\otimes\mathcal G\ar[ul]_{Ad_g\otimes id}\ar[dl]^{\gamma_g}\\
g_!f_!f^\ast(g^!\mathcal F\otimes g^\ast\mathcal
G)\{m\}\ar[r]^{\hspace{0,7cm}g_!Tr_f}
\ar[u]^{g'^!f^\ast\varphi_g}& g_!(g^!\mathcal F\otimes
g^\ast\mathcal G)\ar[u]_{g_!\varphi_g}}
$$

\medskip Here the diagrams (1) and (3) commute by the definitions of
$\varphi$ and $t$, respectively, and the identification
$g'^!=f^!g^!$ is just the one for which the diagrams (2) commute
(these morphisms being defined as the adjoints of $g'_!$, $f_!$
and $g_!$). The commutativity of (4) is easily checked (cf. also
\cite{SGA 4.3} XVII 5.2.4), and (5) commutes by our definition of
(5.7.3), together with an obvious ''associativity'' for
$\varphi_f$. The remaining squares commute by functoriality.

\medskip Now let $Y=Spec~(k)$. Then by definition the map product for
$X$ is the composition
$$
\xymatrix{ H^{-a}(X, a^!_X\mathbb Z/n(-b)\ar@{=}[d]& \otimes &
H^i_Z(X,\,a_X^\ast\mathbb Z/n(j))\ar[r]\ar@{=}[d] & H^{i-a}_Z(X,\,
a^!_X\mathbb Z/n(j-b))\ar@{=}[d]\\ H_a(X,\, \mathbb Z/n(b)) & &
H^i_Z(X,\,\mathbb Z /n(j)) & H_{a-i} (Z,\,\mathbb Z/n(b-j))}
$$
induced by the usual cup product (with supports) and
$\varphi_{a_X}$, together with the identification $H^\cdot_Z(X,\,
a_X^! \mathcal F)=H^\cdot(Z,\, Ri^!a^!_X\mathcal F)= H^\cdot(Z,\,
a_z^!\mathcal F)$ for $i:Z\hookrightarrow X$. Thus 3.3 $(n)$
follows from the commutativity of (5.7.4).

\medskip Next we consider 5.3 (o). Let $i:Z\hookrightarrow X$ be a
closed immersion with open complement $j:U=X-Z\hookrightarrow X$.
For a flat morphism $f:X'\rightarrow X$, which is equidimensional
of dimension $m$, consider the cartesian squares
$$
\xymatrix{ Z'\ar@{^{(}->}[r]^{i'} \ar[d]^{f_Z} & X'\ar[d]^f &
U'\ar@{_{(}->}[l]_{j'}\ar[d]^{f_U}\\
Z\ar@{^{(}->}[r]^{i'} & X & U\ar@{_{(}->}[l]_{j}& .}
$$
The relative sequence 3.1 (e) for the triple $(Z,\, X,\, U)$ is
obtained by taking the cohomology on $X$ of the canonical exact
triangle
$$
i_\ast Ri^! a_X^!\mathbb Z/n(-b)\stackrel{Ad_i}{\longrightarrow}
a_X^!\mathbb Z /n (-b)\stackrel{ad_i}{\longrightarrow} Rj_\ast
j^\ast a^!_X\mathbb Z/n(-b)\rightarrow\; ,
$$
and identifying $Ri^!a^!_X=a^!_Z$ and $j^\ast a_X^!= a^!_U$ (Here
we used that $Ri_\ast =i_\ast $, and have written $Ri^!$ since
$i^!$ may be misinterpreted as the functor '' sections with
support in $Z$'' whose derivative $Ri^!$ is). Hence 5.3 (o).
follows from the fact that one has natural identifications of
exact triangles
$$
\xymatrix{ i'_\ast Ri'^!a_{X'}^!\mathcal F\ar[r]\ar@{=}[d] &
a_{X'}^!\mathcal F\ar[r]\ar@{=}[d] & Rj'_\ast j'^\ast
a_{X'}^!\mathcal F\ar[r]\ar@{=}[d]  & \\
i'_\ast f_Z^! Ri^! a_{X}^!\mathcal F\ar[r]\ar[d]_{\beta '}^{\wr}&
f^! a_{X}^!\mathcal F\ar[r]\ar@{=}[d] & Rj'_\ast f^!_U j^\ast
a_{X}^!\mathcal F\ar[r]\ar[d]_{\beta ''}^{\wr} & \\
f^! i_\ast Ri' a_X^!\mathcal F\ar[r] & f^! a^!_X\mathcal F\ar[r] &
f^! Rj_\ast j^\ast a^!_X\mathcal F\ar[r] & & , }
$$
for any complex of sheaves $\mathcal F$ on $Spec~(k)$, where
$\beta '$ and $\beta ''$ are obtained from the ''adjoint base
change isomorphism'' \cite{SGA 4.3} XVIII 3.1.12.3. Note that by
definition the diagram
$$
\xymatrix{ i'_\ast Ri^! f^! \mathcal G\ar[r]^{Ad_i'}\ar@{=}[d] &
f^!\mathcal G\\i'_\ast f_z^! Ri^!\mathcal G\ar[r]^{\beta '} & f^!
i_\ast Ri^! \mathcal G\ar[u]_{f^!Ad_i}}
$$
commutes, similarly for $\beta ''$ (where $j^!=j^\ast$ and $j'^!
=j'^\ast$).

\medskip Finally we show that the homotopy invariance 5.5 (i) holds for
\'{e}tale homology (this implies property 5.5 (ii) as well, as is
clear from the proof of 5.5). Since $p:\mathbb A^1_X\rightarrow X$
is acyclic (cf. [Mi] VI 4.20), the restriction map
$$
H^{-a}(X,\, a^!_X\mathbb Z/n(-b))\rightarrow H^{-a}(\mathbb
A^1_X,\, p^\ast a_X^!\mathbb Z/n(-b))
$$
is an isomorphism for all $a,\,  b,\,\in\mathbb Z$. We conclude by
recalling that
$$
t_p: p^\ast \mathcal L(1)[2]\rightarrow p^!\mathcal L
$$
is an isomorphism for any $\mathcal L$ in $D_c^b(X,\, \mathbb
Z/n)$ (for this it suffices that $p$ is smooth of relative
dimension 1, cf. \cite{SGA 4.3} XVIII proof of 3.2.5).

\noi We can now collect the fruits of our efforts.

\begin{theo}\label{Theorem5.8} Let $k$ be a field, and let $n$ be a natural
number invertible in $k$. For $i,\, j\in \mathbb Z$ let $\mathcal H^i_n(j)$
be the Zariski sheaf on the category
$\mbox{Sch}^{noeth}/k$ of all noetherian
$k$-schemes associated to the presheaf
$$
U \shortmid\hspace{-0,7mm}\rightsquigarrow H^i(U,\, \mathbb
Z/n(j))
$$
given by \'{e}tale cohomology. Then for all $\nu,\, i,\, j\in
\mathbb Z$ the functor
$$
X \shortmid\hspace{-0,7mm}\rightsquigarrow H^\nu(X,\, \mathcal H^i_n(j))
$$
is a sufficiently rigid functor on $\mbox{Sch}^{noeth}/k$.
\end{theo}

\bigskip\noi {\bf Proof} By Proposition 2.1 and Lemma 3.3, \'{e}tale cohomology is
a sufficiently rigid functor on $\mbox{Sch}^{\scriptsize noeth}/k$. In
view of Theorem 4.1, we then have to show 4.1 (i) and (ii) for the
$\mathcal H^i_n(j)$. Since it suffices to consider the bigger
category $\mbox{Sch}^{\mbox{\scriptsize noeth}}/k_0$, where $k_0$ is the
prime field, we may assume that $k$ is perfect, and by Lemma 4.4 we may
restrict our attention to algebraic $k$-schemes. By Proposition 5.7,
$X\mapsto H^\ast(X,\mathbb Z/n(\ast))$ is part of an
extended Poincar\'{e} duality theory with homotopy invariance.
Therefore the claim follows from Propositions 5.2 and 5.4.

\medskip By Theorem 5.8, we obtain case (7) of Theorem 0.3 in the
introduction.

\begin{coro}\label{Corollary 5.9}
Let $L$ be a complete discrete valuation field, let $n$ be a natural number invertible in $L$,
and let $L^0$ be subfield which is algebraically closed in $L$ and dense in $L$ for the valuation
topology. Then the following holds, with the notations from the proof of Proposition 5.2:

\v\noi (a) If $X$ is any variety over $L^0$, then the induced map
$$
H^\nu(X,\cH^i_n(j)) \rightarrow H^\nu(X\times_{L^0}L,\cH^i_n(j))
$$
is an isomorphism for all $i$ and $j$.

\v\noi (b) If $X$ is any variety over $L^0$, then the induced map
$$
E^2_{d-\nu,a-d}(X) \rightarrow E^2_{d-\nu,a-d}(X\times_{L^0}L)
$$
is an isomorphism for all $\nu$, $a$ and $d$.
\end{coro}

\medskip\noindent \textbf{Proof} Property (a) follows since $X \mapsto H^\nu(X,\cH^i_n(j))$
is a sufficiently rigid functor, and property (b) follows from the isomorphism of long exact sequences
described in (5.27), by induction on dimension. In fact, as in the proof of Theorem 5.8 we may assume that the ground field is
perfect, and then $X$ has an open affine dense subvariety $U$ which is smooth, so that the claim is true
for $U$ by the natural commutative diagram
$$
\xymatrix{H^\nu(U,\cH_a(b)) \ar[r]\ar[d]^{\alpha} & H^\nu(U\times_{L^0}L,\cH_a(b))\ar[d]^{\alpha}\\
          E^2_{d-\nu,a-d}(U)\ar[r]                & E^2_{d-\nu,a-d}(U\times_{L^0}L)}
$$
with the vertical isomorphisms $\alpha$ from (5.2.6). Then we use induction on dimension and the
following commutative diagram with exact rows established in Remark 5.3
$$
\xymatrix{ \ldots\ar[r] & H^\nu_Z(X,\cH_a(b)) \ar[r]\ar[d]^{\alpha} & H^\nu(X,\cH_a(b))\ar[r]\ar[d]^{\alpha} & H^\nu(U,\cH_a(b))
\ar[r]\ar[d]^{\alpha} & H^{\nu+1}_Z(X,\cH_{a-1}(b)) \ar[d]^{\alpha} \ldots\\
\ldots\ar[r] & E^2_{d-\nu,a-d}(Z)\ar[r] & E^2_{d-\nu,a-d}(X)\ar[r] & E^2_{d-\nu,a-d}(U)\ar[r] & E^2_{d-\nu-1,a-d}(Z,b) \ldots }
$$
In fact, this diagram maps to the similar diagram obtained after base change from $L^0$ to $L$, and since these
maps induce isomorphisms at the places for $X$ and $U$, the morphisms give isomorphisms everywhere.

\begin{rem}\label{Remark5.10}
It is shown in \cite{JSS} that the complex $E^1_{d-\nu,a-d}(X)$ can be identified with the complex $C^{r,s}_n(X)$
defined by K. Kato in \cite{Ka} for  $(r,s) = (?,?)$ after changing the differentials by some signs. Therefore the morphism
$$
C^{r,s}_n(X) \rightarrow C^{r,s}_n(X\times_{L^0}L)
$$
is a quasiisomorphism.
\end{rem}

\vspace{0.5cm}

\noi
Uwe Jannsen\\
Fakult\"at f\"ur Mathematik\\
Universit\"at Regensburg\\
93040 Regensburg\\
GERMANY\\
uwe.jannsen@mathematik.uni-regensburg.de

\end{document}